\documentclass[10pt]{article}

\usepackage{amsmath}
\usepackage{amsfonts}
\usepackage{amssymb}
\usepackage{amsthm}
\usepackage{graphicx}
\usepackage{verbatim}
\usepackage{subfig}
\usepackage{caption}
\usepackage{listings}
\usepackage{color}
\usepackage{algorithm}
\usepackage{algpseudocode}
\usepackage{pifont}
\usepackage{hyperref}

\ProvideTextCommand{\DJ}{OT1}{\raisebox{0.25ex}{-}\kern-0.4em D}

\hypersetup{
    colorlinks=true,
    linkcolor=blue,
    filecolor=magenta,      
    urlcolor=cyan,
    pdftitle={Overleaf Example},
    pdfpagemode=FullScreen,
    }

\numberwithin{equation}{section}
\newtheorem{remark}{Remark}
\numberwithin{remark}{section}
\newtheorem{lemma}{Lemma}
\numberwithin{lemma}{section}
\newtheorem{theorem}{Theorem}
\newtheorem{corollary}{Corollary}
\numberwithin{theorem}{section}
\numberwithin{corollary}{section}
\numberwithin{figure}{section}
\numberwithin{table}{section}
\newtheorem{example}{Example}
\numberwithin{example}{section}
\newtheorem{proposition}{Proposition}
\numberwithin{proposition}{section}

\begin{document}
\title{An Unconditionally Energy-Stable and Orthonormality-Preserving Iterative Scheme for the Kohn-Sham Gradient Flow Based Model}

\author{Xiuping Wang\thanks{Computational Transport Phenomena Laboratory (CTPL), Division of Computer, Electrical and Mathematical Sciences and Engineering (CEMSE), King Abdullah University of Science and Technology (KAUST), Saudi Arabia. (Email: xiuping.wang@kaust.edu.sa)}
\and Huangxin Chen\thanks{School of Mathematical Sciences and Fujian Provincial Key Laboratory on Mathematical Modeling and HighPerformance Scientific Computing, Xiamen University, Xiamen, Fujian, People's Republic of China. (Email: chx@xmu.edu.cn)} 
\and Jisheng Kou\thanks{Key Laboratory of Rock Mechanics and Geohazards of Zhejiang Province, Shaoxing University, Shaoxing, People's Republic of China. (Email: jishengkou@163.com)}
  \thanks{School of Mathematics and Statistics, Hubei Engineering University, Xiaogan, Hubei, People's Republic of China.}
 \and Shuyu Sun\thanks{Corresponding author. Computational Transport Phenomena Laboratory (CTPL), Division of Physical Sciences and Engineering (PSE), King Abdullah University of Science and Technology (KAUST), Saudi Arabia. (Email: shuyu.sun@kaust.edu.sa)} 
}

\date{}

\maketitle

\begin{abstract}
	We propose an unconditionally energy-stable, orthonormality-preserving, component-wise splitting iterative scheme for the Kohn-Sham gradient flow based model in the electronic structure calculation. 
	We first study the scheme discretized in time but still continuous in space. The component-wise splitting iterative scheme changes one wave function at a time, similar to the Gauss-Seidel iteration for solving a linear equation system. At the time step $n$, the orthogonality of the wave function being updated to other wave functions is preserved by projecting the gradient of the Kohn-Sham energy onto the subspace orthogonal to all other wave functions known at the current time, while the normalization of this wave function is preserved by projecting the gradient of the Kohn-Sham energy onto the subspace orthogonal to this wave function at $t_{n+1/2}$.  The unconditional energy stability is nontrivial, and it comes from a subtle treatment of the two-electron integral as well as a consistent treatment of the two projections.  Rigorous mathematical derivations are presented to show our proposed scheme indeed satisfies the desired properties. We then study the fully-discretized scheme, where the space is further approximated by a conforming finite element subspace. For the fully-discretized scheme, not only the preservation of orthogonality and normalization (together we called orthonormalization) can be quickly shown using the same idea as for the semi-discretized scheme, but also the highlight property of the scheme, i.e., the unconditional energy stability can be rigorously proven. 
	The scheme allows us to use large time step sizes and deal with small systems involving only a single wave function during each iteration step.
	Several numerical experiments are performed to verify the theoretical analysis, where the number of iterations is indeed greatly reduced as compared to similar examples solved by the Kohn-Sham gradient flow based model in the literature. 
\end{abstract}

\smallskip
\textbf{Keywords}: Density functional theory; gradient flow; orthonormality-preserving schemes; energy stability.

\section{Introduction}\label{Sect. 1}
The Kohn-Sham density functional theory (DFT) is the most widely used approach for many-electron systems in quantum chemistry and physics \cite{kohn1965self}. 
The straightforward way to obtain the lowest-energy solutions is to solve the Kohn-Sham equations \cite{dai2008three,dai2011finite, hu2018multilevel, hu2023accelerating}, a nonlinear eigenvalue problem, where the Hartree-Fock self-consistent field method is commonly applied \cite{ehrenreich1959self, slater1953generalized, matsen2006self, saad2010numerical}. The self-consistent field method assumes that the solution can be represented by a linear combination of atomic orbitals in the entire $\mathbb{R}^3$ space, which may not be accurate enough. It can also be computationally expensive, especially in evaluating the two-electron repulsion integrals for large molecules \cite{barca2021faster, gill_prism_1991, rak_brush_2015}.

\smallskip
Lately, a growing focus has been on the direct energy minimization model. Rather than directly solving the Kohn-Sham equations, the model searches the ground state by minimizing the energy under the Pauli exclusion principle or, equivalently, the orthonormality constraint. 
It is worth mentioning that numerous optimization techniques have been put forth for the Kohn-Sham energy minimization model \cite{dai2017conjugate,francisco2004globally,host2008ground, ulbrich2015proximal, vecharynski2015projected, wen2013feasible,gao2019parallelizable}.
Under appropriate assumptions, the optimization techniques discussed herein yield local convergence.
Nevertheless, satisfying the orthonormality constraint is a challenging task that frequently demands implementing specific techniques; sometimes, they bring significant computational costs to the optimization process.

\smallskip
This paper utilizes the Kohn-Sham gradient flow based model proposed in \cite{dai_gradient_2020}, where an extended gradient is introduced. One notable benefit in contrast to the direct energy minimization approach is that this particular gradient flow model ensures energy stability and orthonormality constraint at the continuous level as well as at the fully discretized level. 
It is imperative to acknowledge the pre-existing literature on employing gradient flow based models in tackling eigenvalue problems \cite{schneider2009direct, bao2004computing, chen1998global}. 

\smallskip
To our knowledge, the existing numerical schemes for the Kohn-Sham gradient flow based model do not possess unconditional energy stability and preserve the orthonormality constraint simultaneously without modifying the original energy \cite{wang2023sav}. Moreover, the existing schemes for the gradient flow model typically arrive at a coupled system, thereby causing a substantial rise in computational costs, particularly as the number of orbitals rises.

\smallskip
Unconditionally energy stable schemes have been actively studied in recent two decades, particularly for multi-component and multi-phase fluid flow and transport problems \cite{feng2022fully, li2017numerical, feng2023energy, shen2022energy, smejkal2022multi, sun2019darcy, sun2020reservoir, li2022fully, fan2020unconditionally, fan2017componentwise, qiao2014two}. Due to the strong nonlinearity, the strong coupling between various physics, and the non-convex feature of the background energy in these problems, small time steps typically have to be used due to stability concerns instead of maintaining accuracy in temporal discretization if the scheme is only conditionally energy stable. For scenarios like these, it is crucial to design unconditionally energy-stable schemes, which not only substantially increase the computational efficiency of the simulation, but also greatly enhance the robustness of the calculation. Unconditionally energy-stable schemes have also been heavily studied for gradient flow problems \cite{yang2017efficient, gong2020arbitrarily, zhang2020unconditionally, shen2019new}; however, to the best of our knowledge, unconditionally energy-stable but orthonormality-preserving schemes for the Kohn-Sham gradient flow-based model in the electronic structure calculation have not been reported yet in the literature.

\smallskip
The contribution of this paper is that we propose an iterative, unconditionally energy-stable, orthonormality-preserving numerical scheme for the Kohn-Sham gradient flow based model. The proposed scheme employs the technique of component-wise splitting, leading to an algorithm that resembles the Gauss-Siedel method for solving systems of linear equations. The properties mentioned above guarantee that we can use large time step sizes in the simulation and solve small systems in each time step, saving massive computational costs compared to other algorithms.  Unlike the SAV-based method \cite{wang2023sav}, our proposal method does not modify the original energy and does not have a tuned parameter for the algorithm. 

\smallskip
The rest of the paper is organized as follows. In the next Section \ref{Sect. 2}, we present the problem setup and the gradient flow formulation at the PDE level.  In Section \ref{Sect. 3}, we propose a new orthonormality-preserving and unconditionally energy-stable time-marching algorithm to solve the gradient flow problem, which yields an orthonormality-preserving, unconditionally energy-stable, and component-splitting iterative scheme for the original energy minimization problem. The rigorous proof is given in this section.  Section \ref{Sect. 4} establishes a fully discrete scheme via the finite element method. In Section \ref{Sect. 5}, some numerical examples are provided to verify the theory analysis. Finally, we give some concluding remarks in Section \ref{Sect. 6}.

\section{Preliminaries}\label{Sect. 2}
The total energy of a system with $N$ orbitals consists of the following parts:
\begin{equation}
	E=E_{\text {kinetic }}+E_{\text {ext }}+E_{\mathrm{Har}}+E_{\mathrm{xc}}+E_{\mathrm{nuc}},
\end{equation}
where $E_{\text {kinetic }}$ is the Kohn-Sham kinetic energy, which is expressed in terms of the orbital $\Psi:=\left\{\psi_1, \psi_2, \cdots, \psi_N\right\}$ as
\begin{equation}
	E_{\text {kinetic }}=\frac{1}{2} \sum_{l=1}^N \int_{\mathbb{R}^3}f_l\left|\nabla \psi_l\right|^2 d \mathbf{r},
\end{equation}
where $f_l$ is the occupation number of the $l$-th orbital. 

\smallskip
$E_{\text {ext }}$ is the external energy expressed as a functional of the external potential $V_{\text{ext}}(\mathbf{r})$ and density $\rho(\mathbf{r})$:
\begin{equation}
	E_{\mathrm{ext}}=\int_{\mathbb{R}^3} V_{\mathrm{ext}}(\mathbf{r}) \rho(\mathbf{r}) d \mathbf{r}, \, V_{\mathrm{ext}}(\mathbf{r})=-\sum_{j=1}^M \frac{Z_j}{\left|\mathbf{r}-\mathbf{R}_j\right|}, \,\rho(\mathbf{r})=\sum_{l=1}^N f_l\psi_l^2(\mathbf{r})
\end{equation} 
where $M$ denotes the number of nuclei, $Z_j$ and $\mathbf{R}_j$ are the corresponding nuclei charge and position of the  $j-$th nuclei, respectively. 
For the simplicity in notations, we present our algorithm below without explicitly listing $f_l$, i.e., we consider the cases that $f_l = \mathrm{constant}, l = 1, ..., N$; in particular, all numerical examples shown in this paper are cases with $f_l = 2, l = 1, ..., N$. 

\smallskip
$E_{\mathrm{Har}}$ is the Hartree energy concerning with the Hartree potential $V_{\mathrm{Har}}$:
\begin{equation}
	E_{\mathrm{Har}}=\frac{1}{2} \int_{\mathbb{R}^3} V_{\mathrm{Har}}(\mathbf{r}) \rho(\mathbf{r}) d \mathbf{r}, \, V_{\mathrm{Har}}([\rho] ; \mathbf{r})=\int_{\mathbb{R}^3} \frac{\rho\left(\mathbf{r}^{\prime}\right)}{\left|\mathbf{r}-\mathbf{r}^{\prime}\right|} d \mathbf{r}^{\prime}
\end{equation}

\smallskip
$E_{\mathrm{xc}}$ is the exchange–correlation energy for which we do not possess the exact formula. It can be expressed as the inner product of the exchange-correlation energy per unit density $\epsilon_{\mathrm{xc}}(\rho)$ and density.
\begin{equation}
	E_{\mathrm{xc}}=\int_{\mathbb{R}^3} \epsilon_{\mathrm{xc}}(\rho) \rho(\mathbf{r}) d \mathbf{r}.
\end{equation}
The last term $E_{\mathrm{nuc}}$ is a constant for the energy of nucleon-nucleon interactions. 
Since the main purpose of this paper is to present and demonstrate our new algorithm, for simplicity, we ignore the last two energy $E_{\mathrm{xc}}$ and $E_{\mathrm{nuc}}$ in this paper. Energy-stable numerical treatment of $E_{\mathrm{xc}}$ will be investigated in our ongoing work. 

\smallskip
The total energy now reads
\begin{equation}
	E=E_{\text {kinetic }}+E_{\text {ext }}+E_{\mathrm{Har}},
\end{equation}
and the solution of ground states can be obtained by minimizing the total energy with orthonormality constraints, which can be formulated as follows:
$$
\begin{aligned}
& \min _{\Psi} E(\Psi) \\
& \text { s.t. }\langle\Psi, \Psi\rangle=I_N,
\end{aligned}
$$
where $\langle\cdot, \cdot\rangle$ stands for the standard $L^2$ inner product, and $I_N$ denotes the $N \times N$ identity matrix. In this paper, we call $\psi_i$ and $\psi_j$ ($i\ne j$) orthogonal to each other if and only if $\langle\psi_i, \psi_j\rangle=0$; we call $\psi_i$ is normalized if and only if $\langle\psi_i, \psi_i\rangle=1$.  Orthogonality and normalization together is called orthonormality. 

\smallskip
For brevity in notations, we first define the operators $H, H_L$, and $H_{\mathrm{Har}}$ by
\begin{equation}
	H \psi_k:=-\frac{1}{2} \nabla^2 \psi_k+\left(V_{\mathrm{ext}}(\mathbf{r})+V_{\mathrm{Har}}([\rho] ; \mathbf{r})\right) \psi_k=H_L \psi_k+H_{\mathrm{Har}}(\rho) \psi_k, \, 1 \le k \le N,
\end{equation}
where $H_L \psi_k:=-\frac{1}{2} \nabla^2 \psi_k+V_{\text {ext }}(\mathbf{r}) \psi_k$ is the linear part, and $H_{\mathrm{Har}}(\rho) \psi_k:=$ $V_{\mathrm{Har}}([\rho] ; \mathbf{r}) \psi_k$ is the nonlinear part.
Then the total energy can be expressed as
\begin{equation}
E=\sum_{l=1}^N\left\langle H_L \psi_l, \psi_l\right\rangle +\frac{1}{2} \sum_{l=1}^N\left\langle H_{\mathrm{Har}}(\rho) \psi_l, \psi_l\right\rangle =\sum_{l=1}^N\left\langle\left(H_L+\frac{1}{2} H_{\mathrm{Har}}\right) \psi_l, \psi_l\right\rangle,	
\end{equation}
and its gradient is
\begin{equation}
	\left(\nabla E(\Psi)\right)_k=2 H \psi_k=2 H_L \psi_k+2 H_{\mathrm{Har}}(\rho) \psi_k.
\end{equation}

We then introduce the Stiefel manifold, defined as follows
\begin{equation} 
	\mathcal{M}^N=\left\{\Psi \in\left(\mathrm{H}^1\left(\mathbb{R}^3\right)\right)^N: \langle\Psi, \Psi\rangle=I_N \right\}.
\end{equation}
The gradient on the Stiefel manifold $\mathcal{M}^N$ of $\left(\nabla E(\Psi)\right)$ at $\Psi$ (see \cite{edelman1998geometry}), denoted by $\nabla_G E(\Psi)$, can be written as
\begin{equation}
	\left(\nabla_G E(\Psi)\right)_k=\left(\nabla E(\Psi)\right)_k-\sum_{l=1}^N\left\langle\left(\nabla E(\Psi)\right)_l, \psi_k\right\rangle \psi_l.
\end{equation}
We note that the gradient on the Stiefel manifold $\mathcal{M}^N$ of $\left(\nabla E(\Psi)\right)$ is the projection of the original gradient $\nabla E(\Psi)$ onto the subspace that is orthogonal to all $\psi_k,k=1,\cdots,N$, because both $H_L$ and $H_{\mathrm{Har}}$ are symmetric operators. 

To propose a gradient flow model that automatically preserves the orthonormal constraint, we have to extend the gradient $\left(\nabla_G E(\Psi)\right)$ from $\mathcal{M}^N$ to $\left(\mathrm{H}^1\left(\mathbb{R}^3\right)\right)^N$. The extended gradient, still denoted as $\left(\nabla_G E(\Psi)\right)$, is defined as
\begin{equation}
	\left(\nabla_G E(\Psi)\right)_k=\sum_{l=1}^N\left\langle\psi_l, \psi_k\right\rangle \left(\nabla E(\Psi)\right)_l-\sum_{l=1}^N\left\langle\left(\nabla E(\Psi)\right)_l, \psi_k\right\rangle \psi_l.
\end{equation}
The corresponding Kohn-Sham gradient flow model can now be formulated as
\begin{equation}\label{Grad_Flow_Model}
	\left\{\begin{array}{l}
\frac{\partial \Psi(\mathbf{r}, \tau)}{\partial \tau}=-\nabla_G E(\Psi), \quad 0< \tau <\infty, \\
\Psi(\mathbf{r}, 0)=\Psi_0,
\end{array}\right.
\end{equation}
where $\Psi_0 = \{ \psi_{0,1},\psi_{0,2},...\psi_{0,N} \}\in \mathcal{M}^N$ is the initial condition. Note that the first equation in (\ref{Grad_Flow_Model}) can also be written as
\begin{equation}
\frac{\partial \psi_k(\mathbf{r}, \tau)}{\partial \tau}=-2 \sum_{l=1}^N\left\langle\psi_l, \psi_k\right\rangle H \psi_l+2 \sum_{l=1}^N\left\langle H \psi_l, \psi_k\right\rangle \psi_l.
\end{equation}
We can do a simple variable substitution of $t = 2\tau$ to obtain
\begin{equation}
\frac{\partial \psi_k(\mathbf{r}, t)}{\partial t}=- \sum_{l=1}^N\left\langle\psi_l, \psi_k\right\rangle H \psi_l+ \sum_{l=1}^N\left\langle H \psi_l, \psi_k\right\rangle \psi_l.
\end{equation}
A direct verification can show that the solution $\Psi$ sits in the Stiefel manifold $\mathcal{M}^N$. Indeed, for any $1 \le m \le N$, we have
\begin{align}
& \frac{d}{d t}\left\langle\psi_k, \psi_m\right\rangle  \notag \\
& =\left\langle\frac{\partial \psi_k}{\partial t}, \psi_m\right\rangle +\left\langle\frac{\partial \psi_m}{\partial t}, \psi_k\right\rangle  \\
& =\left\langle\sum_{l=1}^N\left\langle H \psi_l, \psi_k\right\rangle  \psi_l-\sum_{l=1}^N\left\langle\psi_l, \psi_k\right\rangle  H \psi_l, \psi_m\right\rangle \notag \\
& +\left\langle\sum_{l=1}^N\left\langle H \psi_l, \psi_m\right\rangle  \psi_l-\sum_{l=1}^N\left\langle\psi_l, \psi_m\right\rangle  H \psi_l, \psi_k\right\rangle  \\
& = 0.
\end{align}

\begin{proposition}
	(see Proposition 3.2 in \cite{dai_gradient_2020})The model equations (\ref{Grad_Flow_Model}) satisfy
	\begin{enumerate}
		\item $\langle\Psi, \Psi\rangle=I_N$;
		\item $\frac{dE}{dt} \le 0$.
	\end{enumerate}
	
\end{proposition}

\section{An Iterative Orthonormality-Preserving and Unconditionally Energy-Stable Algorithm}\label{Sect. 3}
This section introduces an iterative scheme that preserves orthonormality and unconditional discrete energy stability.  Let 
\begin{equation}
	\left\{t_n: n=0,1,2 \cdots\right\} \subset[0,+\infty)
\end{equation}
  be discrete points such that $0=t_0<t_1<t_2<\cdots<t_n<\cdots$ and set $\Delta t_n=t_{n+1}-t_n$. For any $n \ge 0$ and $1 \le k \le N$, we introduce the following modified midpoint scheme:
\begin{equation}\label{Num_Scheme}
	\frac{\psi_k^{n+1}-\psi_k^n}{ \Delta t_n} = \mathcal{A}(\psi_k^{n+\frac{1}{2}}), 
\end{equation}
where we have
\begin{align}
	\psi_k^{n+\frac{1}{2}} &=\frac{1}{2}\left(\psi_k^n+\psi_k^{n+1}\right), \\
	\mathcal{A}(\psi_k^{n+\frac{1}{2}}) &= -\left\langle\psi_k^{n+\frac{1}{2}}, \psi_k^{n+\frac{1}{2}}\right\rangle   H^{n+\frac{1}{2}} \psi_k^{n+\frac{1}{2}}+\left\langle H^{n+\frac{1}{2}} \psi_k^{n+\frac{1}{2}}, \psi_k^{n+\frac{1}{2}}\right\rangle   \psi_k^{n+\frac{1}{2}} \notag \\
	&+ \sum_{l \neq k}\left\langle\psi_k^{n+\frac{1}{2}}, \psi_k^{n+\frac{1}{2}}\right\rangle  \left\langle H^{n+\frac{1}{2}} \psi_k^{n+\frac{1}{2}}, \psi_l^{n+\delta_l}\right\rangle   \psi_l^{n+\delta_l} \\
	H^{n+\frac{1}{2}} &= H\left(\psi_1^{n+1}, \cdots, \psi_{k-1}^{n+1}, \psi_k^{n+1}, \psi_{k+1}^n, \cdots, \psi_N^n\right) \\
	\delta_l &= \left\{\begin{array}{l}
		1, \, l < k \\
		0, \, l > k
	\end{array} \right. 
\end{align}

We note that we called the scheme ``modified midpoint scheme'' because in a component-wise midpoint scheme, the term $H^{n+\frac{1}{2}}$ above should read  
$$H^{n+\frac{1}{2}} = H\left(\psi_1^{n+1}, \cdots, \psi_{k-1}^{n+1}, \psi_k^{n+1/2}, \psi_{k+1}^n, \cdots, \psi_N^n\right).$$
However, even though this component-wise midpoint scheme also preserves orthogonality and normalization, it does not have unconditional energy stability, as we will see from the proof later in this paper.  We also note that notation-wise, we should call it $H_k^{n+\frac{1}{2}}$ as it also depends on $k$, but since this dependence is clear from the context, we abuse the notation a little and omit $k$ for simplicity. 

Let us denote $\Psi^{n+1} := \left\{\psi_1^{n+1}, \psi_2^{n+1}, \cdots, \psi_N^{n+1}\right\}$.
We present the following two theorems to show the scheme satisfies the desired properties.

\begin{theorem}\label{Semi_orthonormal}
	The numerical solution of the scheme (\ref{Num_Scheme}) satisfies
	\begin{equation}
		\langle\Psi^{n+1}, \Psi^{n+1}\rangle=I_N, \, \forall n \ge 0,
	\end{equation}
	if we have $\langle\Psi_0, \Psi_0\rangle=I_N$.
	\begin{proof}
		We first take the inner product of the scheme (\ref{Num_Scheme}) with $\psi_j^{n+\delta_j}(j\neq k)$.
		\begin{align}\label{psi_k psi_l property}
			\begin{aligned}
			\left\langle \frac{\psi_k^{n+1}-\psi_k^n}{ \Delta t_n},\psi_j^{n+\delta_j} \right\rangle = \left\langle \frac{\psi_k^{n+1}}{ \Delta t_n},\psi_j^{n+\delta_j} \right\rangle = \left\langle\mathcal{A}(\psi_k^{n+\frac{1}{2}})), \psi_j^{n+\delta_j} \right\rangle.
			\end{aligned}
		\end{align}
		Then we can obtain
		\begin{align}
			\begin{aligned}
			\left\langle \mathcal{A}(\psi_k^{n+\frac{1}{2}}), \psi_j^{n+\delta_j} \right\rangle &=  -\left\langle\psi_k^{n+\frac{1}{2}}, \psi_k^{n+\frac{1}{2}}\right\rangle  \left\langle H^{n+\frac{1}{2}} \psi_k^{n+\frac{1}{2}},  \psi_j^{n+\delta_j}\right\rangle  \\
			&+ \left\langle H^{n+\frac{1}{2}} \psi_k^{n+\frac{1}{2}}, \psi_k^{n+\frac{1}{2}}\right\rangle   \left\langle \psi_k^{n+\frac{1}{2}}, \psi_j^{n+\delta_j} \right\rangle  \\
			&+ \sum_{l \neq k}\left\langle\psi_k^{n+\frac{1}{2}}, \psi_k^{n+\frac{1}{2}}\right\rangle  \left\langle H^{n+\frac{1}{2}} \psi_k^{n+\frac{1}{2}}, \psi_l^{n+\delta_l}\right\rangle   \left\langle \psi_l^{n+\delta_l}, \psi_j^{n+\delta_j}\right\rangle.
			\end{aligned}
		\end{align}
		Note that we have $\left\langle \psi_l^{n+\delta_l}, \psi_j^{n+\delta_j}\right\rangle = \delta_{lj}$, then we can arrive at
		\begin{align}\label{A_property}
			\begin{aligned}
			\left\langle \mathcal{A}(\psi_k^{n+\frac{1}{2}}), \psi_j^{n+\delta_j} \right\rangle &=  -\left\langle\psi_k^{n+\frac{1}{2}}, \psi_k^{n+\frac{1}{2}}\right\rangle  \left\langle H^{n+\frac{1}{2}} \psi_k^{n+\frac{1}{2}},  \psi_j^{n+\delta_j}\right\rangle  \\
			&+ \left\langle H^{n+\frac{1}{2}} \psi_k^{n+\frac{1}{2}}, \psi_k^{n+\frac{1}{2}}\right\rangle   \left\langle \psi_k^{n+\frac{1}{2}}, \psi_j^{n+\delta_j} \right\rangle  \\
			&+\left\langle\psi_k^{n+\frac{1}{2}}, \psi_k^{n+\frac{1}{2}}\right\rangle  \left\langle H^{n+\frac{1}{2}} \psi_k^{n+\frac{1}{2}}, \psi_j^{n+\delta_j}\right\rangle    \\
			&= \left\langle H^{n+\frac{1}{2}} \psi_k^{n+\frac{1}{2}}, \psi_k^{n+\frac{1}{2}}\right\rangle   \left\langle \psi_k^{n+\frac{1}{2}}, \psi_j^{n+\delta_j} \right\rangle  \\
			&= \frac{1}{2}  \left\langle H^{n+\frac{1}{2}} \psi_k^{n+\frac{1}{2}}, \psi_k^{n+\frac{1}{2}}\right\rangle   \left\langle \psi_k^{n+1}, \psi_j^{n+\delta_j} \right\rangle
			\end{aligned}
		\end{align}
		Substituting the equation (\ref{A_property}) into (\ref{psi_k psi_l property}), we can obtain
		\begin{equation}
			\left\langle \psi_k^{n+1}, \psi_j^{n+\delta_j} \right\rangle = 0,
		\end{equation}		
		unless 
		\begin{equation}\label{continuity condition}
			\frac{1}{ \Delta t_n}=\frac{1}{2}\left\langle H^{n+\frac{1}{2}} \psi_k^{n+\frac{1}{2}}, \psi_k^{n+\frac{1}{2}}\right\rangle.
		\end{equation}
		By continuity argument, we can still show $\left\langle \psi_k^{n+1}, \psi_j^{n+\delta_j} \right\rangle = 0$ even if (\ref{continuity condition}) holds. To show its normalization property, we now take the inner product of the scheme (\ref{Num_Scheme}) again with $\psi_k^{n+\frac{1}{2}}$ to obtain
		\begin{align}\label{psi_k preserving property}
			\begin{aligned}
				\left\langle \frac{\psi_k^{n+1}-\psi_k^n}{ \Delta t_n}, \psi_k^{n+\frac{1}{2}} \right\rangle = \frac{1}{2\Delta t_n } \left(\left\langle \psi_k^{n+1},\psi_k^{n+1} \right\rangle - 1\right)= \left\langle\mathcal{A}(\psi_k^{n+\frac{1}{2}}), \psi_k^{n+\frac{1}{2}} \right\rangle,
			\end{aligned}
		\end{align}
		where we have
		\begin{align}\label{A preserving}
			\begin{aligned}
				\left\langle \mathcal{A}(\psi_k^{n+\frac{1}{2}}), \psi_k^{n+\frac{1}{2}} \right\rangle &=  -\left\langle\psi_k^{n+\frac{1}{2}}, \psi_k^{n+\frac{1}{2}}\right\rangle  \left\langle H^{n+\frac{1}{2}} \psi_k^{n+\frac{1}{2}},  \psi_k^{n+\frac{1}{2}}\right\rangle \\
			&+ \left\langle H^{n+\frac{1}{2}} \psi_k^{n+\frac{1}{2}}, \psi_k^{n+\frac{1}{2}}\right\rangle   \left\langle \psi_k^{n+\frac{1}{2}}, \psi_k^{n+\frac{1}{2}} \right\rangle   \\
			&+ \sum_{l \neq k}\left\langle\psi_k^{n+\frac{1}{2}}, \psi_k^{n+\frac{1}{2}}\right\rangle  \left\langle H^{n+\frac{1}{2}} \psi_k^{n+\frac{1}{2}}, \psi_l^{n+\delta_l}\right\rangle   \left\langle \psi_l^{n+\delta_l}, \psi_k^{n+\frac{1}{2}}\right\rangle \\
			&= 0.
			\end{aligned}
		\end{align}
		Combining the results of equations (\ref{psi_k preserving property})-(\ref{A preserving}), we arrive at
		\begin{equation}
			\left\langle \psi_k^{n+1},\psi_k^{n+1} \right\rangle= 1,
		\end{equation}
		which concludes the proof.
	\end{proof}
\end{theorem}
Before we state the next theorem, we recall the expression of the total energy and rewrite it as
\begin{align}\label{EnergyDef}
	\begin{aligned}
	& E(\Psi)=E\left(\left\{\psi_1, \cdots, \psi_{k-1}, \psi_k, \psi_{k+1}, \cdots, \psi_N\right\}\right)  \\
= & \sum_{l=1}^N\left\langle\left(H_L+\frac{1}{2} H_{\mathrm{Har}}(\Psi)\right) \psi_l, \psi_l\right\rangle \\
= & \sum_{l=1}^N\left\langle H_L \psi_l, \psi_l\right\rangle + \frac{1}{2} \sum_{l=1}^N \sum_{m=1}^N\left\langle L_{\mathrm{Har}} \psi_l^2, \psi_m^2\right\rangle,
	\end{aligned}
\end{align}
where we utilize the definition of density and define
\begin{equation}
	\left\langle L_{\mathrm{Har}} \rho_l, \rho_m\right\rangle := \int_{\mathbb{R}^3} \int_{\mathbb{R}^3} \frac{\rho_l(\mathbf{r}) \rho_m\left(\mathbf{r}^{\prime}\right)}{\left|\mathbf{r}-\mathbf{r}^{\prime}\right|} d \mathbf{r} d \mathbf{r}^{\prime}.
\end{equation}
\begin{lemma}\label{L_property}
	The following two inequalities hold
	\begin{align}
		\left\langle L_{\mathrm{Har}} \rho_2,\left(\rho_1-\rho_2\right)\right\rangle &\leq\left\langle L_{\mathrm{Har}} \rho_1,\left(\rho_1-\rho_2\right)\right\rangle, \\
		\left\langle L_{\mathrm{Har}} \rho_2, \rho_2\right\rangle-\left\langle L_{\mathrm{Har}} \rho_1, \rho_1\right\rangle &\leq 2\left\langle L_{\mathrm{Har}} \rho_2,\left(\rho_2-\rho_1\right)\right\rangle.\label{L_property2}
	\end{align}
	\begin{proof}
		Note that $\left\langle L_{\mathrm{Har}} \rho_l, \rho_m\right\rangle$ is in a symmetric and positive-definite form. Thus, we have
		\begin{align}
			\begin{aligned}
				0 & \leq\left\langle L_{\mathrm{Har}}\left(\rho_1-\rho_2\right),\left(\rho_1-\rho_2\right)\right\rangle \\
				& =\left\langle L_{\mathrm{Har}} \rho_1,\left(\rho_1-\rho_2\right)\right\rangle-\left\langle L_{\mathrm{Har}} \rho_2,\left(\rho_1-\rho_2\right)\right\rangle, 
			\end{aligned}
		\end{align}
		and 
		\begin{equation}
			\begin{aligned}
				& \left\langle L_{\mathrm{Har}} \rho_2, \rho_2\right\rangle-\left\langle L_{\mathrm{Har}} \rho_1, \rho_1\right\rangle \\
				= & 2\left\langle L_{\mathrm{Har}} \rho_2,\left(\rho_2-\rho_1\right)\right\rangle-\left\langle L_{\mathrm{Har}}\left(\rho_1-\rho_2\right),\left(\rho_1-\rho_2\right)\right\rangle \\
				\leq & 2\left\langle L_{\mathrm{Har}} \rho_2,\left(\rho_2-\rho_1\right)\right\rangle .
			\end{aligned}
		\end{equation}
		This completes the proof.
	\end{proof}
\end{lemma}

\begin{lemma}\label{EnergyInequality_lemma}
	For any $n \ge 0$ and $1 \le k \le N$, the following inequality holds
	\begin{equation}
		\frac{E(\Psi^{n+1}_k) - E(\Psi^{n+1}_{k-1})}{2 \Delta t_n} \le \left\langle H^{n+\frac{1}{2}} \psi_k^{n+\frac{1}{2}}, \frac{\psi_k^{n+1}-\psi_k^n}{ \Delta t_n}\right\rangle,
	\end{equation}
	where
	\begin{align}
		\Psi^{n+1}_k &:= \left\{\psi_1^{n+1}, \cdots, \psi_{k-1}^{n+1}, \psi_k^{n+1}, \psi_{k+1}^n, \cdots, \psi_N^n\right\}, \\
		\Psi^{n+1}_{k-1} &:= \left\{\psi_1^{n+1}, \cdots, \psi_{k-1}^{n+1}, \psi_k^{n}, \psi_{k+1}^n, \cdots, \psi_N^n\right\}.
	\end{align}
	\begin{proof}
		By the definition of the total energy (\ref{EnergyDef}), we have
		\begin{align}\label{EnergyInequality}
			\begin{aligned}
& \frac{E(\Psi^{n+1}_k) - E(\Psi^{n+1}_{k-1})}{2( \Delta t_n)} \\
&= \frac{1}{2 \Delta t_n}\left(\sum_{l=1}^k\left\langle H_L \psi_l^{n+1}, \psi_l^{n+1}\right\rangle+\sum_{l=k+1}^N\left\langle H_L \psi_l^{n}, \psi_l^{n}\right\rangle\right) \\
& -\frac{1}{2 \Delta t_n}\left(\sum_{l=1}^{k-1}\left\langle H_L \psi_l^{n+1}, \psi_l^{n+1}\right\rangle+\sum_{l=k}^N\left\langle H_L \psi_l^{n}, \psi_l^{n}\right\rangle\right) \\
& +\frac{1}{4 \Delta t_n}\left\langle L_{\operatorname{Har}}\left(\sum_{l=1}^k\left(\psi_l^{n+1}\right)^2+\sum_{l=k+1}^N\left(\psi_l^n\right)^2\right),\left(\sum_{m=1}^k\left(\psi_m^{n+1}\right)^2+\sum_{m=k+1}^N\left(\psi_m^n\right)^2\right)\right\rangle \\
& -\frac{1}{4 \Delta t_n}\left\langle L_{\operatorname{Har}}\left(\sum_{l=1}^{k-1}\left(\psi_l^{n+1}\right)^2+\sum_{l=k}^N\left(\psi_l^n\right)^2\right),\left(\sum_{m=1}^{k-1}\left(\psi_m^{n+1}\right)^2+\sum_{m=k}^N\left(\psi_m^n\right)^2\right)\right\rangle .
\end{aligned}
		\end{align}
		The first two terms on the right-hand side of (\ref{EnergyInequality}) can be simplified as
		\begin{equation}
			\begin{aligned}
& \frac{1}{2 \Delta t_n}\left(\sum_{l=1}^k\left\langle H_L \psi_l^{n+1}, \psi_l^{n+1}\right\rangle+\sum_{l=k+1}^N\left\langle H_L \psi_l^{n}, \psi_l^{n}\right\rangle\right) \\
& -\frac{1}{2 \Delta t_n}\left(\sum_{l=1}^{k-1}\left\langle H_L \psi_l^{n+1}, \psi_l^{n+1}\right\rangle+\sum_{l=k}^N\left\langle H_L \psi_l^{n}, \psi_l^{n}\right\rangle\right) \\
= & \frac{1}{2 \Delta t_n}\left\langle H_L \psi_k^{n+1}, \psi_k^{n+1}\right\rangle-\frac{1}{2 \Delta t_n}\left\langle H_L \psi_k^n, \psi_k^n\right\rangle \\
= & \left\langle\frac{H_L\left(\psi_k^{n+1}-\psi_k^n\right)}{ \Delta t_n}, \psi_k^{n+\frac{1}{2}}\right\rangle=\left\langle H_L \psi_k^{n+\frac{1}{2}}, \frac{\psi_k^{n+1}-\psi_k^n}{ \Delta t_n},\right\rangle .
\end{aligned}
		\end{equation}
		By using the inequality (\ref{L_property2}) in Lemma \ref{L_property}, we can obtain
		\begin{equation}
			\begin{aligned}
& \frac{1}{4 \Delta t_n}\left\langle L_{\operatorname{Har}}\left(\sum_{l=1}^k\left(\psi_l^{n+1}\right)^2+\sum_{l=k+1}^N\left(\psi_l^n\right)^2\right),\left(\sum_{m=1}^k\left(\psi_m^{n+1}\right)^2+\sum_{m=k+1}^N\left(\psi_m^n\right)^2\right)\right\rangle \\ 
& -\frac{1}{4 \Delta t_n}\left\langle L_{\operatorname{Har}}\left(\sum_{l=1}^{k-1}\left(\psi_l^{n+1}\right)^2+\sum_{l=k}^N\left(\psi_l^n\right)^2\right),\left(\sum_{m=1}^{k-1}\left(\psi_m^{n+1}\right)^2+\sum_{m=k}^N\left(\psi_m^n\right)^2\right)\right\rangle \\
& =\frac{1}{2 \Delta t_n}\left\langle L_{\mathrm{Har}}\left(\sum_{m=1}^{k-1}\left(\psi_m^{n+1}\right)^2+\sum_{m=k+1}^N\left(\psi_m^n\right)^2\right),\left(\psi_k^{n+1}\right)^2-\left(\psi_k^n\right)^2\right\rangle  \\
& +\frac{1}{4 \Delta t_n}\left(\left\langle L_{\mathrm{Har}}\left(\psi_k^{n+1}\right)^2,\left(\psi_k^{n+1}\right)^2\right\rangle -\left\langle L_{\mathrm{Har}}\left(\psi_k^n\right)^2,\left(\psi_k^n\right)^2\right\rangle \right) \\
& \leq \frac{1}{2 \Delta t_n}\left\langle L_{\operatorname{Har}}\left(\sum_{m=1}^{k-1}\left(\psi_m^{n+1}\right)^2+\sum_{m=k+1}^N\left(\psi_m^n\right)^2\right),\left(\psi_k^{n+1}\right)^2-\left(\psi_k^n\right)^2\right\rangle \\
& +\frac{1}{2 \Delta t_n}\left\langle L_{\operatorname{Har}}\left(\psi_k^{n+1}\right)^2,\left(\psi_k^{n+1}\right)^2-\left(\psi_k^n\right)^2\right\rangle  \\
& =\frac{1}{2 \Delta t_n}\left\langle L_{\operatorname{Har}}\left(\sum_{m=1}^k\left(\psi_m^{n+1}\right)^2+\sum_{m=k+1}^N\left(\psi_m^n\right)^2\right),\left(\psi_k^{n+1}\right)^2-\left(\psi_k^n\right)^2\right\rangle \\
& =\frac{1}{ \Delta t_n}\left\langle L_{\operatorname{Har}}\left(\sum_{m=1}^k\left(\psi_m^{n+1}\right)^2+\sum_{m=k+1}^N\left(\psi_m^n\right)^2\right),\left(\psi_k^{n+1}-\psi_k^n\right) \psi_k^{n+\frac{1}{2}}\right\rangle  \\
& =\left\langle H_{\mathrm{Har}}^{n+\frac{1}{2}} \psi_k^{n+\frac{1}{2}}, \frac{\psi_k^{n+1}-\psi_k^n}{ \Delta t_n}\right\rangle, 
\end{aligned}
		\end{equation}
		which concludes the proof.
	\end{proof}
\end{lemma}

\begin{theorem}\label{Unconditionally energy-stable}
	The numerical solution of the scheme (\ref{Num_Scheme}) satisfies
	\begin{equation}
		E(\Psi^{n+1}) \le E(\Psi^{n}), \, \forall n \ge 0,
	\end{equation}
	\begin{proof}
		By utilizing the result in Lemma \ref{EnergyInequality_lemma}, we can arrive at
		\begin{align}
			&\frac{E(\Psi^{n+1}_k) - E(\Psi^{n+1}_{k-1})}{2 \Delta t_n} \le \left\langle H^{n+\frac{1}{2}} \psi_k^{n+\frac{1}{2}}, \frac{\psi_k^{n+1}-\psi_k^n}{ \Delta t_n}\right\rangle \\
			&= \left\langle H^{n+\frac{1}{2}} \psi_k^{n+\frac{1}{2}}, \frac{\psi_k^{n+1}-\psi_k^n}{ \Delta t_n}\right\rangle \\
			&= - \left\langle\psi_k^{n+\frac{1}{2}}, \psi_k^{n+\frac{1}{2}}\right\rangle \Vert H^{n+\frac{1}{2}} \psi_k^{n+\frac{1}{2}} \Vert^2 + \left\langle H^{n+\frac{1}{2}} \psi_k^{n+\frac{1}{2}}, \psi_k^{n+\frac{1}{2}}\right\rangle^2 \notag \\
			&+\left\langle\psi_k^{n+\frac{1}{2}}, \psi_k^{n+\frac{1}{2}}\right\rangle \sum_{l=1}^{k-1}\left\langle H^{n+\frac{1}{2}} \psi_k^{n+\frac{1}{2}}, \psi_l^{n+1}\right\rangle^2 \notag \\
			&+ \left\langle\psi_k^{n+\frac{1}{2}}, \psi_k^{n+\frac{1}{2}}\right\rangle \sum_{l=k + 1}^N\left\langle H^{n+\frac{1}{2}} \psi_k^{n+\frac{1}{2}}, \psi_l^n\right\rangle^2 \\
			&= -c_0^2\left\langle H^{n+\frac{1}{2}} \psi_k^{n+\frac{1}{2}}, H^{n+\frac{1}{2}} \psi_k^{n+\frac{1}{2}}-c_{1, k} \widehat{\psi_k^{n+\frac{1}{2}}}-\sum_{l=1}^{k-1} c_{2, l} \psi_l^{n+1}-\sum_{l=k+1}^N c_{3, l} \psi_l^n\right\rangle \label{Energy_Projection},
		\end{align}
		where we have
		\begin{align}
			c_0^2:&=\left\|\psi_k^{n+\frac{1}{2}}\right\|^2=\frac{1}{2}\left(1+\left\langle\psi_k^n, \psi_k^{n+1}\right\rangle\right) \\
			\widehat{\psi_k^{n+\frac{1}{2}}}&:=\frac{\psi_k^{n+\frac{1}{2}}}{\left\|\psi_k^{n+\frac{1}{2}}\right\|}=\frac{1}{c_0} \psi_k^{n+\frac{1}{2}} \\
			c_{1, k}&:=\left\langle H^{n+\frac{1}{2}} \psi_k^{n+\frac{1}{2}}, \widehat{\psi_k^{n+\frac{1}{2}}}\right\rangle \\
			c_{2, l}&:=\left\langle H^{n+\frac{1}{2}} \psi_k^{n+\frac{1}{2}}, \psi_l^{n+1}\right\rangle \\
			c_{3, l}&:=\left\langle H^{n+\frac{1}{2}} \psi_k^{n+\frac{1}{2}}, \psi_l^n\right\rangle
		\end{align}
		Note that the last three terms in Eq. (\ref{Energy_Projection}) are indeed the projection from $H^{n+\frac{1}{2}} \psi_k^{n+\frac{1}{2}}$ to the space spanned by $\{\psi_1^{n+1},...\psi_{k-1}^{n+1}, \widehat{\psi_k^{n+\frac{1}{2}}}, \psi_{k+1}^n...\psi_N^n\}$. Thus, we can have the following identity
		\begin{align}
			&\left\langle H^{n+\frac{1}{2}} \psi_k^{n+\frac{1}{2}}, H^{n+\frac{1}{2}} \psi_k^{n+\frac{1}{2}}-c_{1, k} \widehat{\psi_k^{n+\frac{1}{2}}}-\sum_{l=1}^{k-1} c_{2, l} \psi_l^{n+1}-\sum_{l=k+1}^N c_{3, l} \psi_l^n\right\rangle \notag \\
			&= \Vert H^{n+\frac{1}{2}} \psi_k^{n+\frac{1}{2}}-c_{1, k} \widehat{\psi_k^{n+\frac{1}{2}}}-\sum_{l=1}^{k-1} c_{2, l} \psi_l^{n+1}-\sum_{l=k+1}^N c_{3, l} \psi_l^n \Vert^2.
		\end{align}
		By substituting the above equation back to Eq. (\ref{Energy_Projection}), we arrive at
		\begin{align}
			&\frac{E(\Psi^{n+1}_k) - E(\Psi^{n+1}_{k-1})}{2 \Delta t_n} \notag \\
			&\le -c_0^2\Vert H^{n+\frac{1}{2}} \psi_k^{n+\frac{1}{2}}-c_{1, k} \widehat{\psi_k^{n+\frac{1}{2}}}-\sum_{l=1}^{k-1} c_{2, l} \psi_l^{n+1}-\sum_{l=k+1}^N c_{3, l} \psi_l^n \Vert^2 \\
			&\le 0,
		\end{align}
		which completes the proof.
	\end{proof}
\end{theorem}
We remark that $\left\|H^{n+\frac{1}{2}} \psi_k^{n+\frac{1}{2}}-c_{1, k} \widehat{\psi_k^{n+\frac{1}{2}}}-\sum_{l=1}^{k-1} c_{2, k} \psi_l^{n+1}-\sum_{l=k+1}^N c_{2, k} \psi_l^n\right\|=$ 0 if and only if the operator $H^{n+\frac{1}{2}}$ maps $\psi_k^{n+\frac{1}{2}}$ onto the current solution space of $\operatorname{span}\left\{\psi_1^{n+1}, \cdots, \psi_{k-1}^{n+1}, \psi_k^{n+\frac{1}{2}}, \psi_{k+1}^n, \cdots, \psi_N^n\right\}$. We further note that when converging, $H^{n+1}$ maps $\operatorname{span}\left\{\psi_1^{n+1}, \cdots, \psi_k^{n+1}, \psi_k^{n+1}, \psi_{k+1}^{n+1}, \cdots, \psi_N^{n+1}\right\}$ onto the same space of $\operatorname{span}\left\{\psi_1^{n+1}, \cdots, \psi_k^{n+1}, \psi_k^{n+1}, \psi_{k+1}^{n+1}, \cdots, \psi_N^{n+1}\right\}$, which implies that $\psi_1^{n+1}, \cdots, \psi_k^{n+1}, \psi_k^{n+1}, \psi_{k+1}^{n+1}, \cdots, \psi_N^{n+1}$ is the solution of the original Kohn-Sham eigenvalue problem unique up to the equivalence relation defined by the orthonormal transform.

\section{Fully Discrete Schemes}\label{Sect. 4}
This section introduces the fully discrete scheme provided by the standard finite element method. Let $\mathcal{T}_h$ be a tetrahedral finite element partition of a finite computational domain $\Omega$ in $\mathbb{R}^3$ with a mesh size $h$. The corresponding finite element space $V_h$ is defined as follows:
\begin{equation}
	V_h=\left\{v\in C^0(\Omega):\left.v\right|_K \in \mathcal{P}_r(K), K \in \mathcal{T}_h\right\},
\end{equation}
where $\mathcal{P}_r$ denotes the space of polynomials of total degree at most $r$.

\smallskip
For simplicity, we consider the homogenous Dirichlet boundary condition.
We let $V^0_h$ be the subspace of $V_h$ satisfying the homogenous Dirichlet boundary condition.  Then we can derive the fully discrete scheme: For $n \ge 0$ and $1 \le k \le N$, find $\psi_{k,h}^{n+1} \in V^0_h$, such that 
\begin{equation}\label{Fully_Num_Scheme}
	\left\langle \frac{\psi_{k,h}^{n+1}-\psi_{k,h}^n}{ \Delta t_n}, v_h \right\rangle = \left\langle \mathcal{A}_h(\psi_{k,h}^{n+\frac{1}{2}}), v_h \right\rangle, \, \forall v_h \in V^0_h,
\end{equation}
where we have
\begin{align}
	\psi_{k,h}^{n+\frac{1}{2}} &=\frac{1}{2}\left(\psi_{k,h}^n+\psi_{k,h}^{n+1}\right), \\
	\mathcal{A}_h(\psi_{k,h}^{n+\frac{1}{2}}) &= -\left\langle\psi_{k,h}^{n+\frac{1}{2}}, \psi_{k,h}^{n+\frac{1}{2}}\right\rangle   H^{n+\frac{1}{2}}_h \psi_{k,h}^{n+\frac{1}{2}}+\left\langle H^{n+\frac{1}{2}}_h \psi_{k,h}^{n+\frac{1}{2}}, \psi_{k,h}^{n+\frac{1}{2}}\right\rangle   \psi_{k,h}^{n+\frac{1}{2}} \notag \\
	&+ \sum_{l \neq k}\left\langle\psi_{k,h}^{n+\frac{1}{2}}, \psi_{k,h}^{n+\frac{1}{2}}\right\rangle  \left\langle H^{n+\frac{1}{2}}_h \psi_{k,h}^{n+\frac{1}{2}}, \psi_{l,h}^{n+\delta_l}\right\rangle   \psi_{l,h}^{n+\delta_l}, \\
	H^{n+\frac{1}{2}}_h &= H_h\left(\psi_{1,h}^{n+1}, \cdots, \psi_{k-1, h}^{n+1}, \psi_{k,h}^{n+1}, \psi_{k+1,h}^n, \cdots, \psi_{N,h}^n\right), \\
	\left\langle H_h \psi_k, v_h \right\rangle &= \left\langle \nabla \psi_k, \nabla v_h \right\rangle+\left\langle \left(V_{\mathrm{ext}}(\mathbf{r})+V_{\mathrm{Har}}([\rho_h] ; \mathbf{r})\right) \psi_k, v_h \right\rangle.
\end{align}
Let us denote $\Psi^{n+1}_h := \left\{\psi_{1,h}^{n+1}, \psi_{2,h}^{n+1}, \cdots, \psi_{N,h}^{n+1}\right\}$. We present the following corollary to show the fully discrete scheme maintains the desired properties.
\begin{corollary}
	For any $n \ge 0$, the numerical solution of the scheme \ref{Fully_Num_Scheme} satisfies
	\begin{equation}
		E(\Psi^{n+1}_h) \le E(\Psi^{n}_h), \, \langle\Psi^{n+1}_h, \Psi^{n+1}_h\rangle=I_N.
	\end{equation}
	\begin{proof}
		The proof can be conducted in the same way as in Theorem \ref{Semi_orthonormal} and Theorem \ref{Unconditionally energy-stable}.
	\end{proof}
\end{corollary}

We also remark that since we are using a conforming Galerkin method for the spatial discretization, i.e., $V^0_h=V^0_h(\Omega)$ is a subspace of $H^1(\mathbb{R}^3)$, we quickly see that when converging, our numerical energy is an upper-bound estimate of the true energy if the exact solution of the PDE approaches to the global minimum solution of the Kohn-Sham energy functional: 
$$\lim_{n\rightarrow\infty}E(\Psi^{n}_h) \ge \lim_{n\rightarrow\infty} E(\Psi^{n}).$$

The energy error $E_{E} := \lim\limits_{n\rightarrow\infty}E(\Psi^{n}_h) - \lim\limits_{n\rightarrow\infty} E(\Psi^{n})$ typically decreases as we refine the mesh.  The analysis of the energy error (and the error on wave functions) is out of the scope of this paper. 

\section{Numerical Implementations}\label{Sect. 5}

This section presents a few numerical examples to validate our theoretical analysis. Note that all the quantities are in the atomic unit. Let the occupation number $f_l = 2$ for all orbitals, and the number of orbitals $N$ satisfies $N = \frac{N_e}{2}$, where $N_e$ is the number of electrons. The number of electrons $N_e$ should be the same as the total nuclei charge for charge-neutral systems.

\smallskip
The stopping criterion is set as $\vert E(\Psi^{n+1}_h) - E(\Psi^{n}_h)\vert \le 1e-6 $. A fixed point iteration with a criterion of $1e-8$ is applied to solve the nonlinear problem. Let $j \ge 0$ be the number of iterations. The detailed scheme can be formulated as
\begin{equation}\label{Ite_Fully_Num_Scheme}
	\left\langle \frac{\psi_{k,h}^{n+1, j+1}-\psi_{k,h}^n}{ \Delta t_n}, v_h \right\rangle = \left\langle \mathcal{A}_h^j(\psi_{k,h}^{n+\frac{1}{2}}), v_h \right\rangle, \, \forall v_h \in V_h^0,
\end{equation}
where we have
\begin{align}
	\psi_{k,h}^{n+\frac{1}{2}, j+1} &=\frac{1}{2}\left(\psi_{k,h}^n+\psi_{k,h}^{n+1, j+1}\right), \\
	\mathcal{A}_h^j(\psi_{k,h}^{n+\frac{1}{2}}) &= -\left\langle\psi_{k,h}^{n+\frac{1}{2},j}, \psi_{k,h}^{n+\frac{1}{2},j}\right\rangle   H^{n+\frac{1}{2},j}_h \psi_{k,h}^{n+\frac{1}{2},j+1} \\
	&+\left\langle H^{n+\frac{1}{2},j}_h \psi_{k,h}^{n+\frac{1}{2}, j}, \psi_{k,h}^{n+\frac{1}{2},j}\right\rangle   \psi_{k,h}^{n+\frac{1}{2},j+1} \notag \\
	&+ \sum_{l \neq k}\left\langle\psi_{k,h}^{n+\frac{1}{2},j}, \psi_{k,h}^{n+\frac{1}{2},j}\right\rangle  \left\langle H^{n+\frac{1}{2},j}_h \psi_{k,h}^{n+\frac{1}{2},j}, \psi_{l,h}^{n+\delta_l}\right\rangle   \psi_{l,h}^{n+\delta_l}, \\
	H^{n+\frac{1}{2},j}_h &= H_h\left(\psi_{1,h}^{n+1}, \cdots, \psi_{k-1, h}^{n+1}, \psi_{k,h}^{n+1, j}, \psi_{k+1,h}^n, \cdots, \psi_{N,h}^n\right), 
\end{align}
where the initial guess $\psi_{k,h}^{n+1, 0} = \psi_{k,h}^{n}$.

\begin{figure}[h!]
	\subfloat{\includegraphics[scale = 0.22]{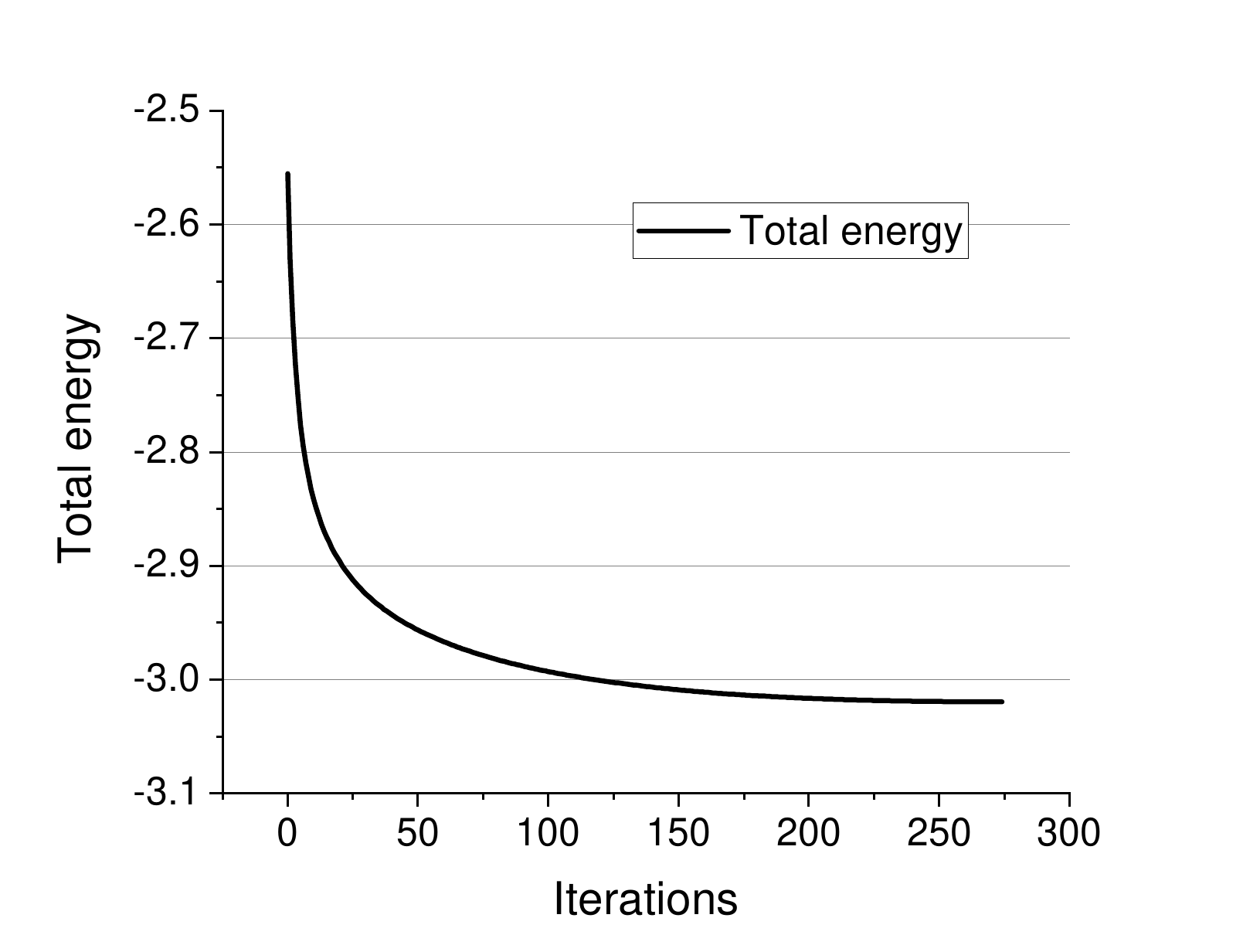}} 
	\hspace{0cm}
	\subfloat{\includegraphics[scale = 0.22]{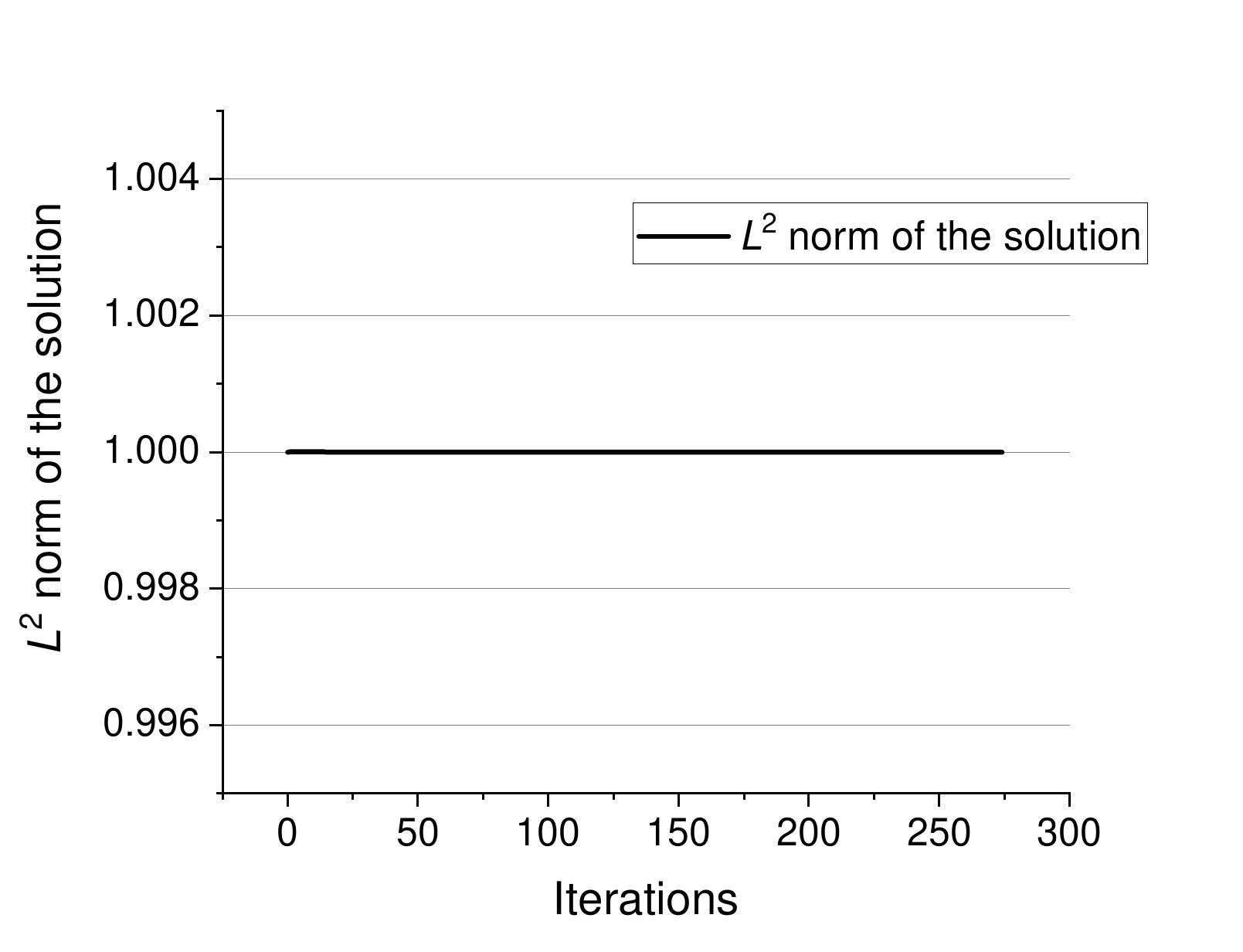}} 
	\caption{The evolution of the computed total energy (in Hartree) with the time step (left) and the $L^2$ norm of the numerical solution (right) for the electronic structure of a helium atom (Example \ref{Ex.1}), demonstrating that our scheme preserves normalization exactly while being strictly energy stable even with large time steps.}
	\label{Ex.1 Fig}
\end{figure}

\begin{figure}[h!]
	\subfloat{\includegraphics[height = 5cm, width = 5cm]{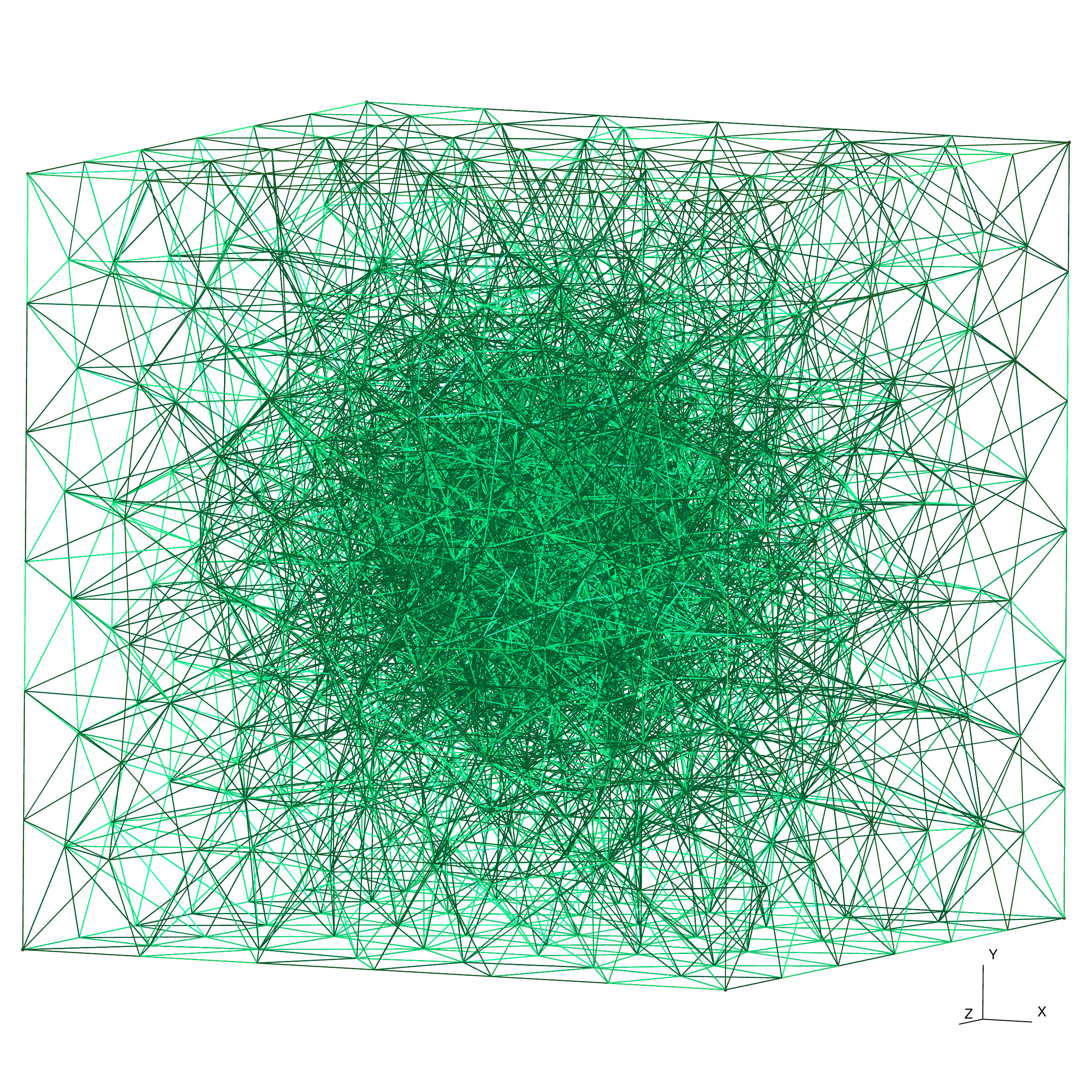}} 
	\hspace{1.5cm}
	\subfloat{\includegraphics[scale = 0.1]{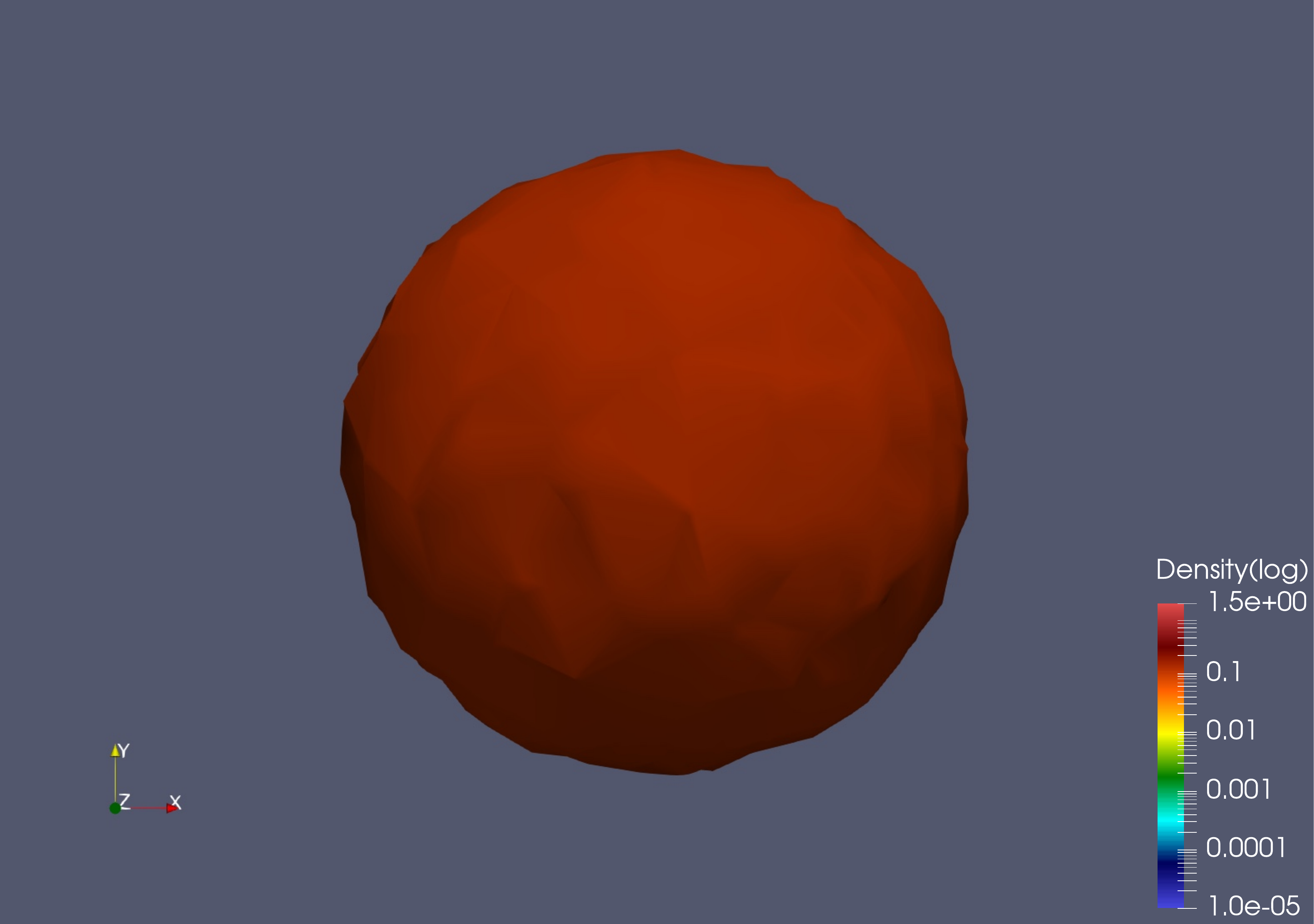}} 
	\caption{The nonuniform mesh used for domain discretization with a total number of degree of freedoms 5400 (left) and the contour plot (right) for the predicted electronic structure of a helium atom (Example \ref{Ex.1}).}
	\label{Ex.1 mesh}
\end{figure}

\begin{figure}
	\subfloat{\includegraphics[scale = 0.1]{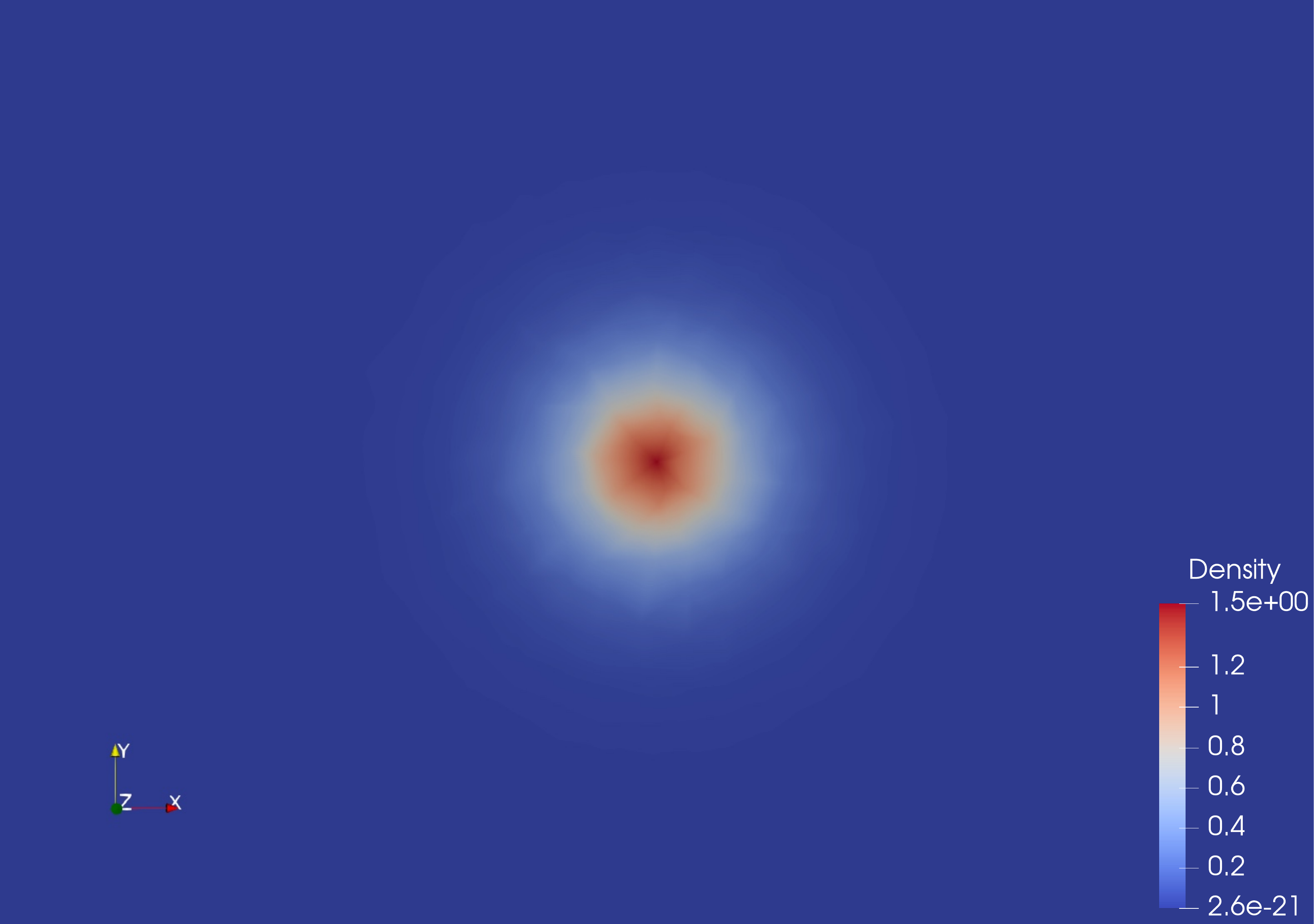}} 
	\hspace{1cm}
	\subfloat{\includegraphics[scale = 0.1]{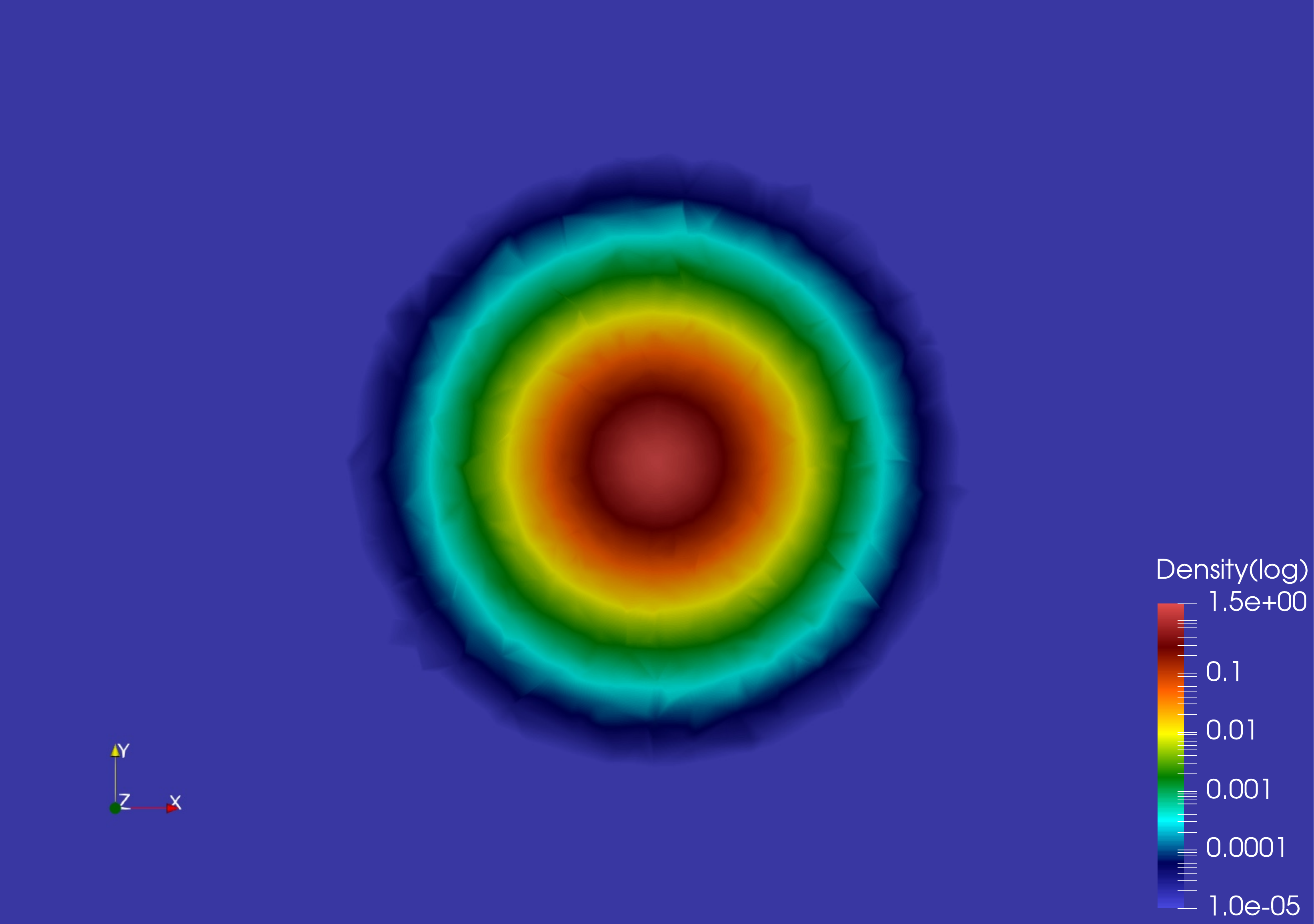}}
	\caption{Profiles of the computed electron density function in a linear scale (left) and in a log scale (right) for the $X Y$ cross-section of the electronic structure of a helium atom (Example \ref{Ex.1})}
	\label{Ex.1 result}
\end{figure}

All the meshes are generated by the Gmsh \cite{geuzaine2009gmsh}. The polynomial order $r$ is fixed as 1.
\begin{example}\label{Ex.1}
\end{example}
We first consider a helium (He) atom with charge two and location $(0,0,0)$. The nonuniform mesh size $h$ is determined by
\begin{equation}
	h(x,y,z) = \frac{x^2+y^2+z^2}{400} + 0.2.
\end{equation}
The initial condition is given as
\begin{equation}
	\psi_{0,1}(\mathbf{r}) = \frac{e^{-2\vert \mathbf{r} - \mathbf{R}_1\vert + 0.5}}{\Vert e^{-2\vert \mathbf{r} - \mathbf{R}_1\vert + 0.5} \Vert}.
\end{equation}
Other parameters are set as follows:
\begin{equation}
	\Omega = [-10, 10]^3, \, \Delta t_n = \Delta t = 1e-4, \, \forall n \ge 0.
\end{equation}
As shown in Figure \ref{Ex.1 Fig}, we conclude that the proposed scheme preserves the discrete energy stability and orthonormal constraint. 
The nonuniform mesh, contour plot, and electron density are shown in Figure \ref{Ex.1 mesh} and Figure \ref{Ex.1 result}, respectively. A spherical shape can be seen from these two figures. 

\begin{figure}[h!]
	\subfloat{\includegraphics[scale = 0.22]{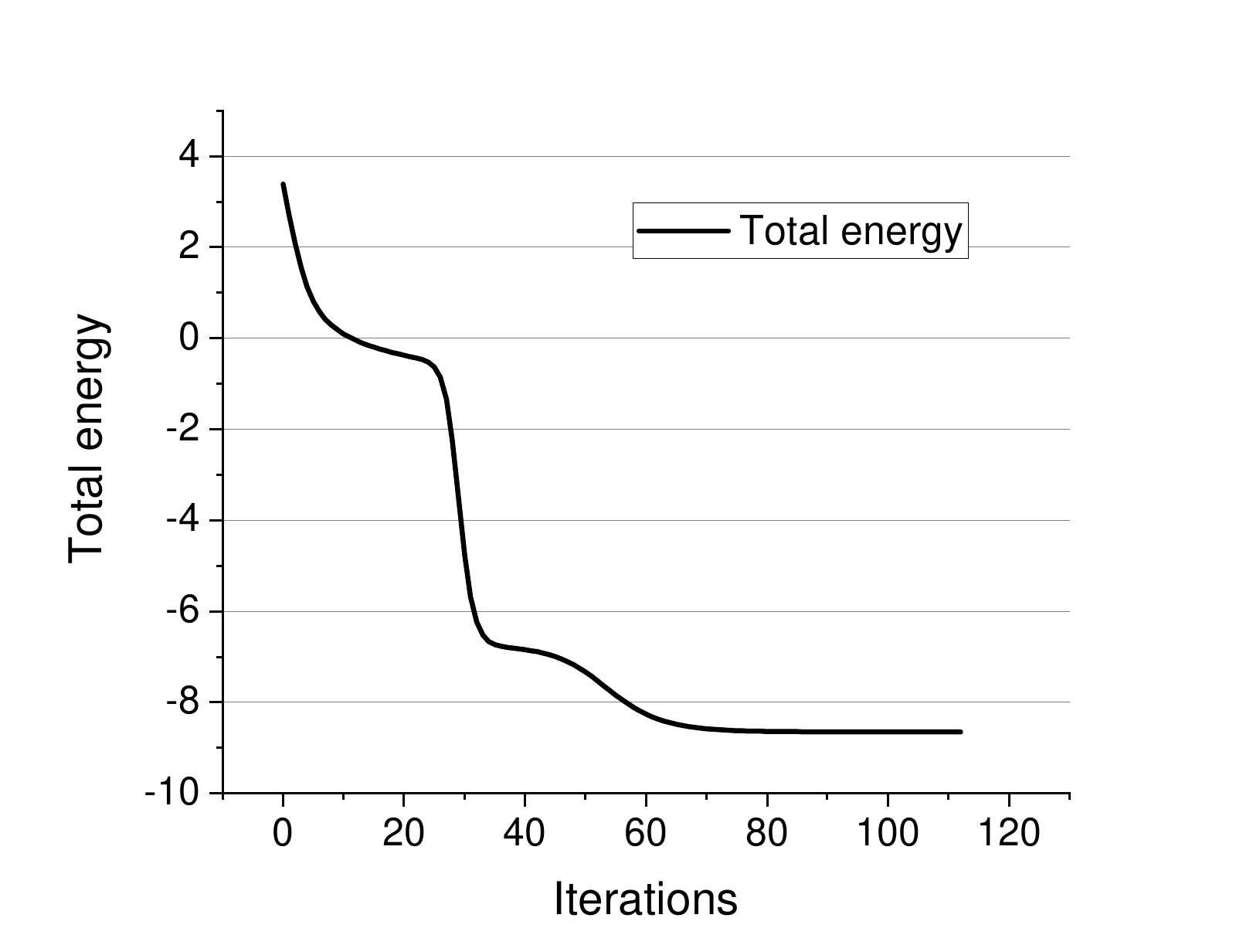}} 
	\hspace{0cm}
	\subfloat{\includegraphics[scale = 0.22]{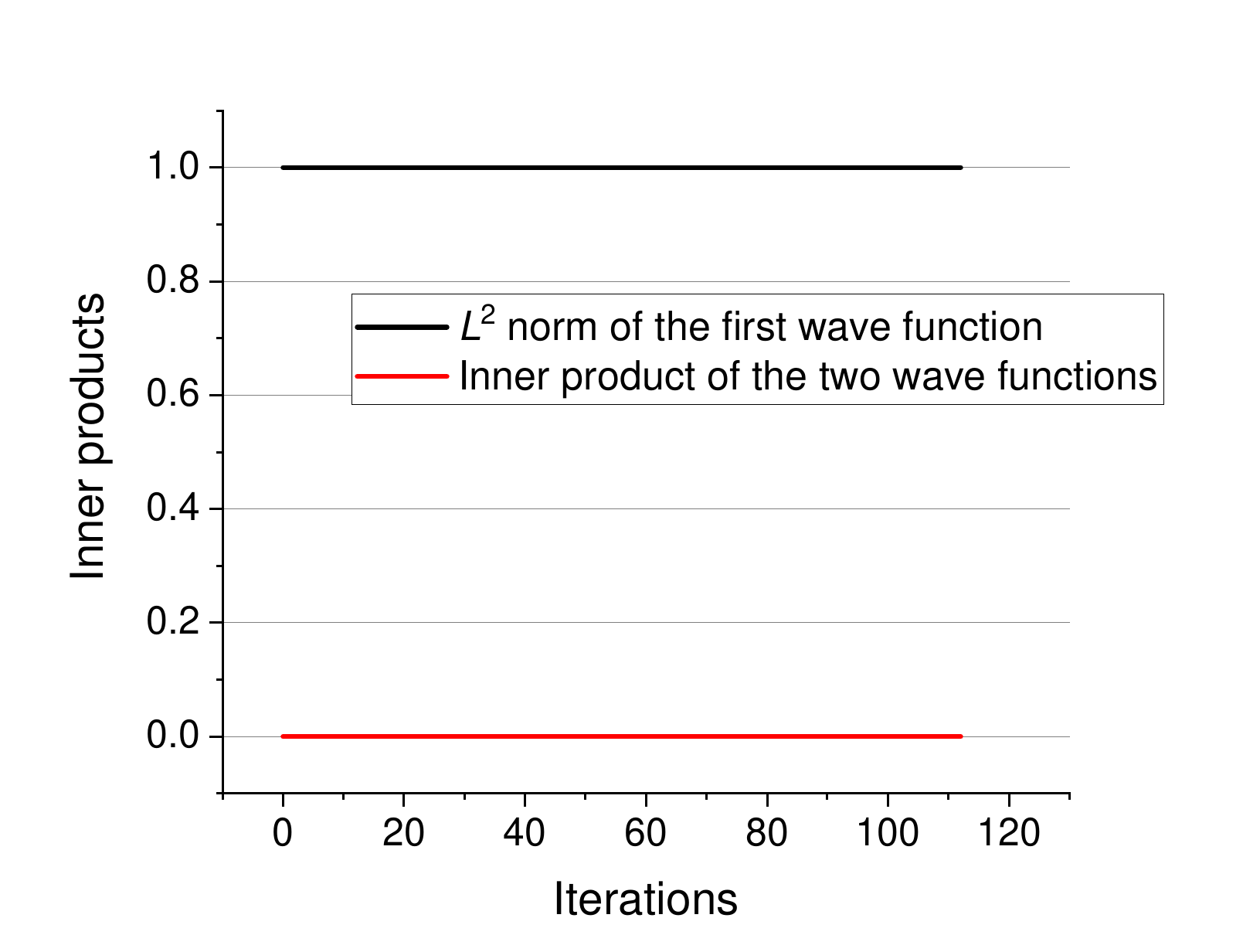}} 
	\caption{The evolution of the computed total energy (in Hartree) with the time step (left) and the $L^2$ norm of the first wave function and inner product of the two wave functions (right) for the electronic structure of a lithium hydride molecule (Example \ref{Ex.2}), demonstrating that our scheme preserves normalization and orthogonality exactly while being strictly energy stable even with large time steps.}
	\label{Ex.2 Fig}
\end{figure}

\begin{figure}[h!]
	\subfloat{\includegraphics[height = 5cm, width = 5cm]{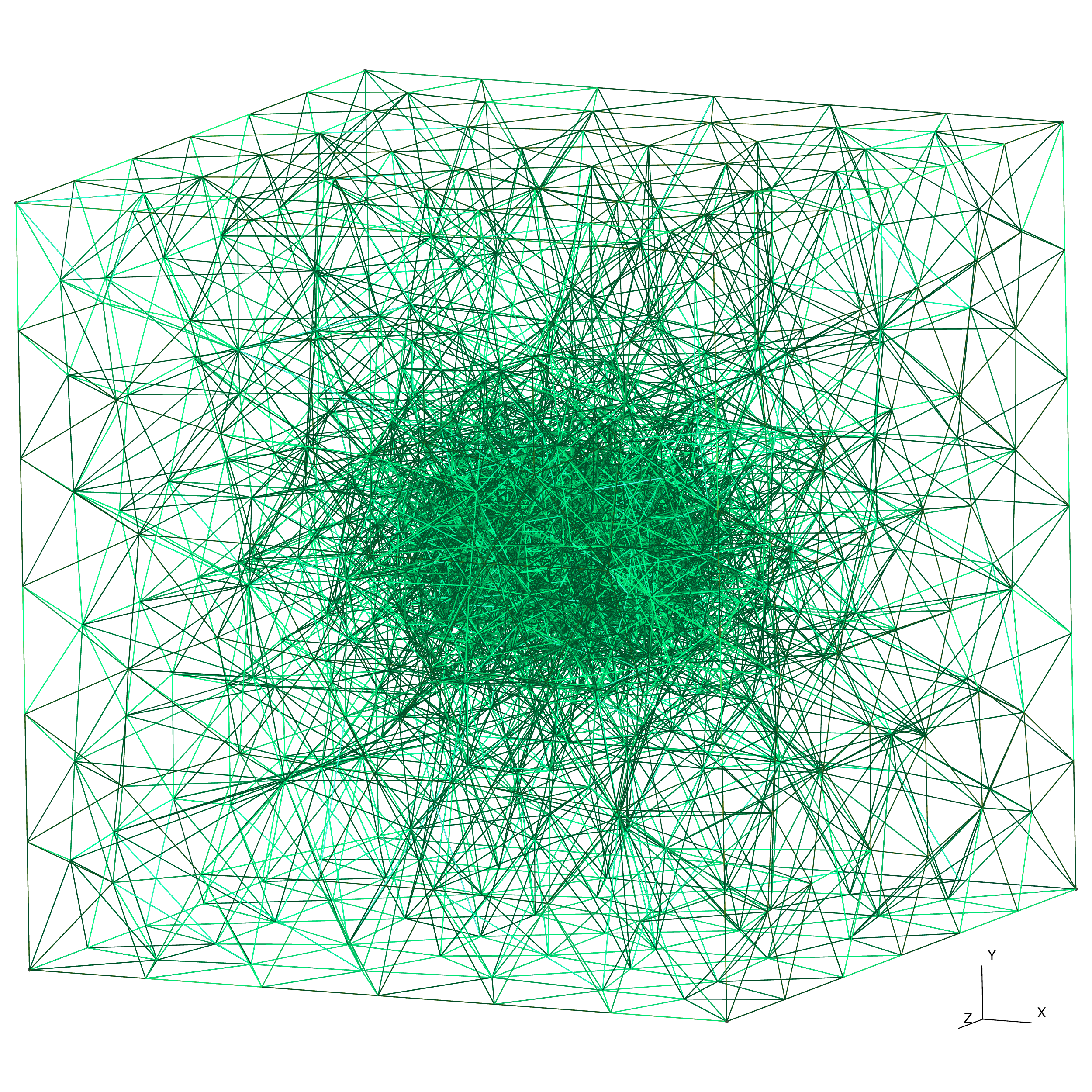}} 
	\hspace{1.5cm}
	\subfloat{\includegraphics[scale = 0.1]{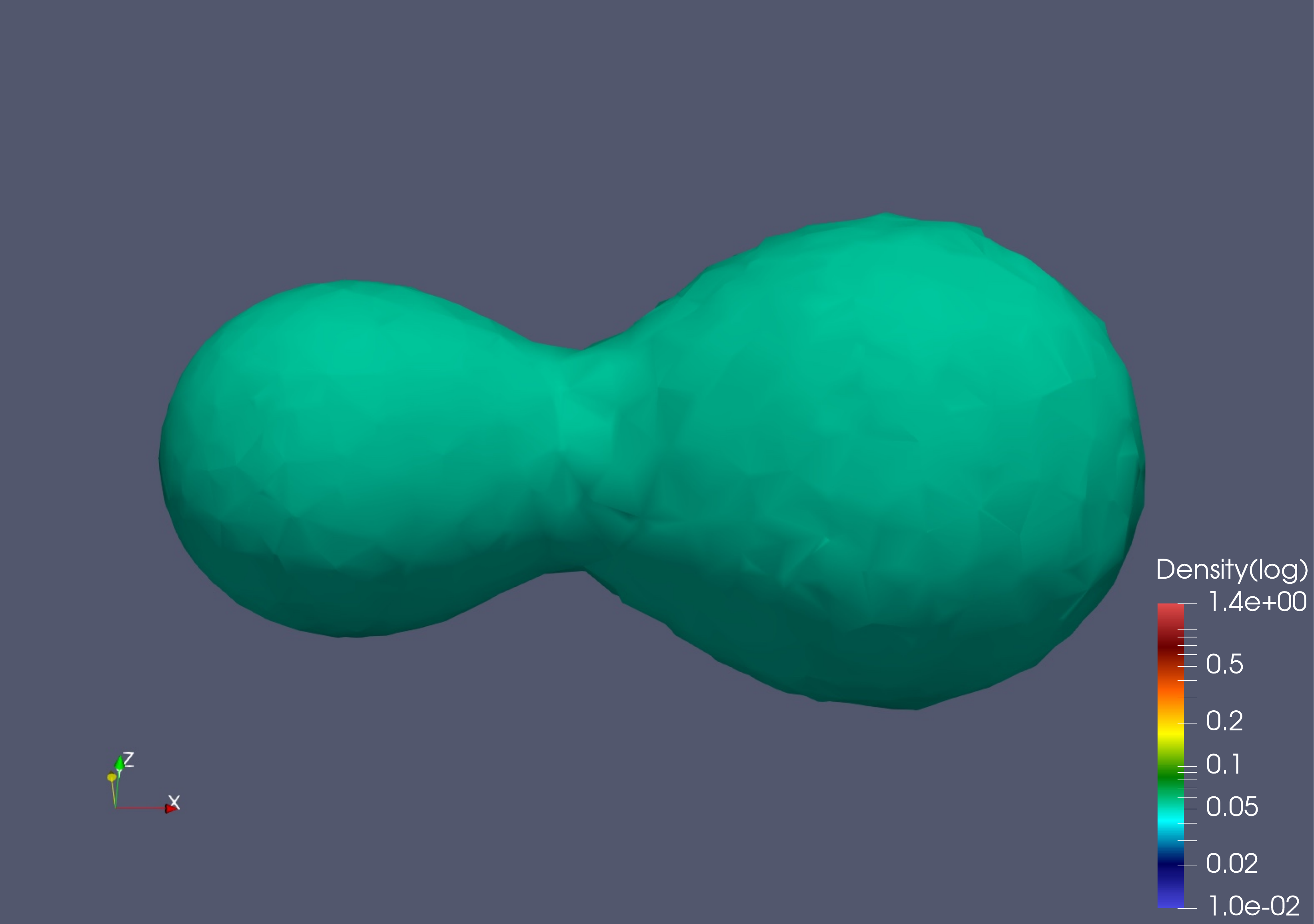}} 
	\caption{The nonuniform mesh used for domain discretization with a total number of degree of freedoms 6909 (left) and the 3D contour plot (right) for the predicted electronic structure of a lithium hydride molecule (Example \ref{Ex.2}).}
	\label{Ex.2 mesh}
\end{figure}

\begin{figure}
	\subfloat{\includegraphics[scale = 0.1]{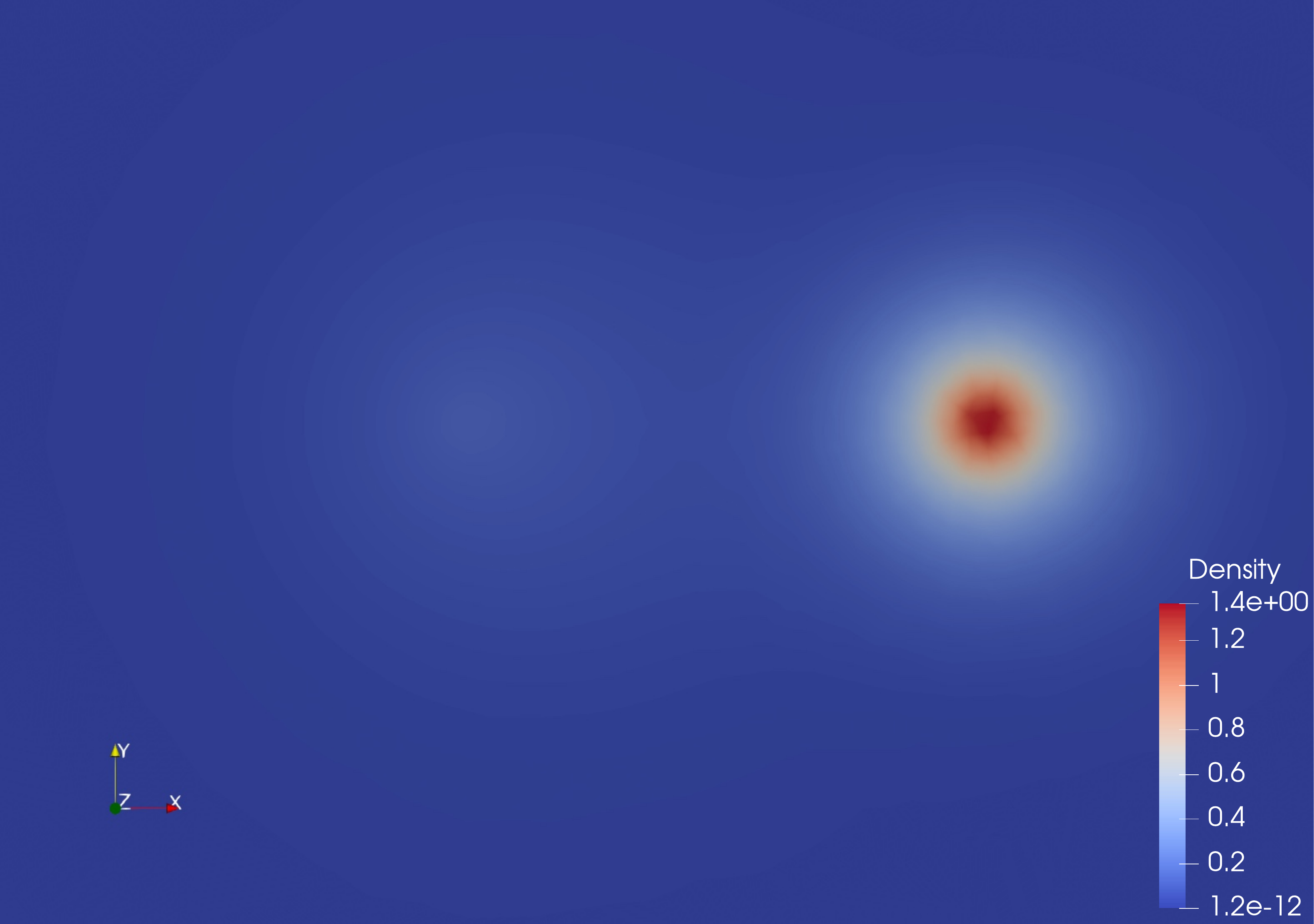}} 
	\hspace{1cm}
	\subfloat{\includegraphics[scale = 0.1]{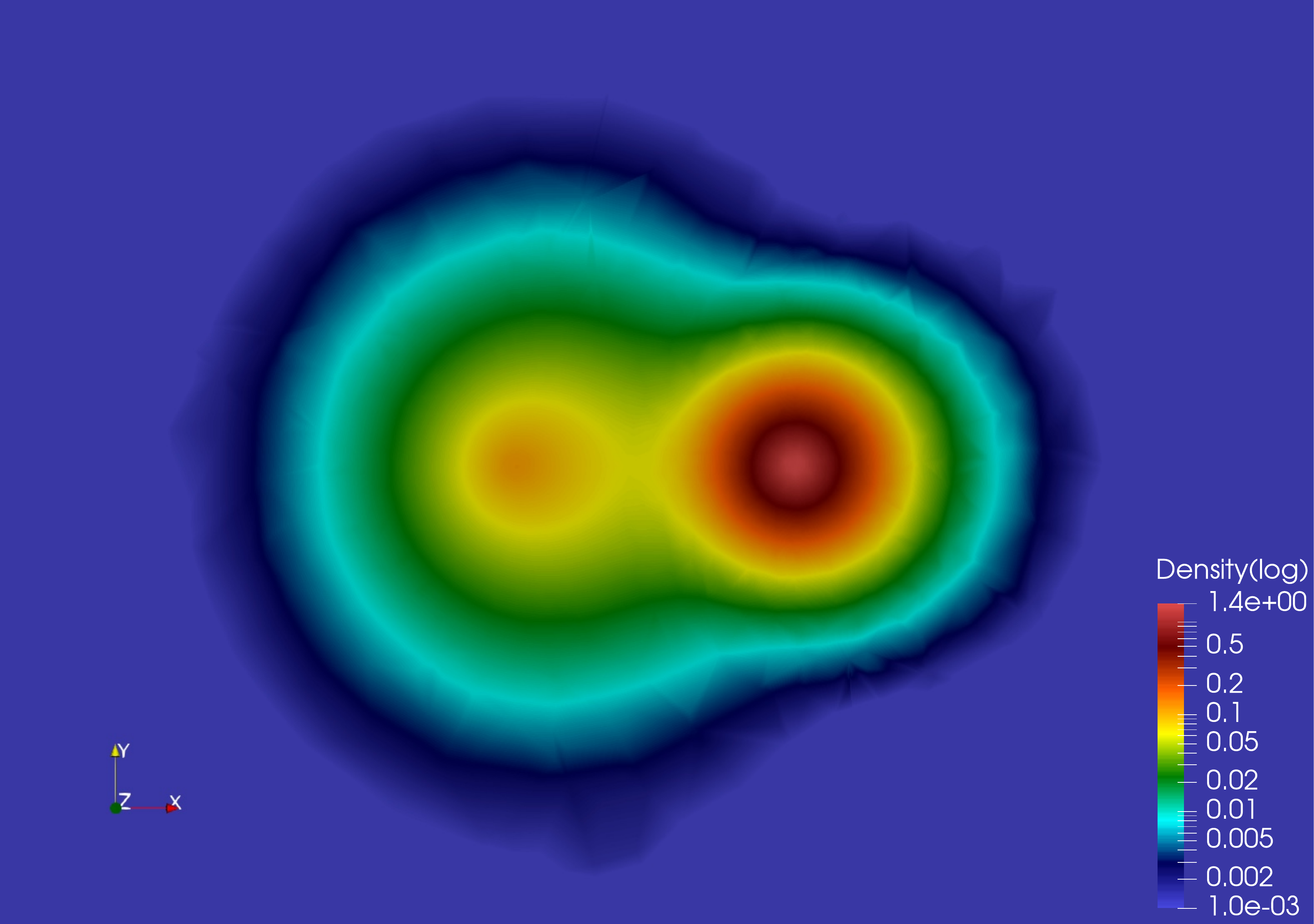}}
	\caption{Profiles of the computed electron density function in a linear scale (left) and in a log scale (right) for the $X Y$ cross-section of the electronic structure of a lithium hydride molecule (Example \ref{Ex.2}).}
	\label{Ex.2 result}
\end{figure}

\begin{example}\label{Ex.2}
\end{example}
In this example, we consider the lithium hydride (LiH) of 2 nuclei with charge $\{1,3\}$ and positions $(-1.0075,0,0)$, $(2.0075,0,0)$. 
The initial condition can be given as
\begin{align}
	\psi_{0,1}(\mathbf{r}) = (1,0,...,0)^T, \,\psi_{0,2}(\mathbf{r}) = (0,1,...,0)^T, \, \psi_{0,l}(\mathbf{r}) \in \mathbb{R}^{N_{\text{dof}} \times 1},
\end{align}
where $N_{\text{dof}}$ is the degrees of freedom. We must ensure that the first and second points are not in the same element. The next step is to normalize the initial condition to let them sit in the Stiefel manifold. This initial condition is designed to be far away from the equilibrium. We utilize a relatively large time step size 
\begin{equation}
	\Delta t_n = \Delta t = 1e-1.
\end{equation}
The computational domain remains the same as the last. The nonuniform mesh size is calculated by
\begin{equation}
	h(x,y,z) = \min(\frac{\min\left((x-1)^2 + y^2 + z^2, (x+2)^2 + y^2 + z^2\right)}{15} + 0.1, h_{\max}),
\end{equation}
where $h_{\max}$ is the default parameter in the Gmsh.
As illustrated in Figure \ref{Ex.2 Fig}, we can still observe the discrete energy stability and orthonormality-preserving property. The simulation only takes 110 steps to converge thanks to the unconditional energy stability, even though the initial condition is quite distant from the ground state. Additionally, it can be observed that the numerical solution consistently resides within the Stiefel manifold, so initiating a restart procedure is unnecessary to ensure orthonormality. We illustrate the nonuniform mesh, profiles of the density, and 3D contour plot in Figure \ref{Ex.2 mesh} and Figure \ref{Ex.2 result}.


\begin{example}\label{Ex.3}
\end{example}

Finally, we simulate a methane ($\text{CH}_4$) molecule of 5 nuclei with charge $\{6,1,1,1,1\}$ and positions $(0,0,0)$, $(c,c,c)$, $(-c,-c,c)$, $(c,-c,-c)$, $(-c,c,-c)$, where the constant $c = 1.1892$\footnote{From Computational Chemistry Comparison and Benchmark DataBase: \url{https://cccbdb.nist.gov/expgeom2.asp?casno=74828&charge=0}.}. The computational domain stays the same as Example \ref{Ex.1}. We give the initial condition in the same way as in the last example, except for the first wave function
\begin{equation}
	\psi_{0,1}(\mathbf{r}) = \frac{e^{-2\vert \mathbf{r} - \mathbf{R}_1\vert}}{\Vert e^{-2\vert \mathbf{r} - \mathbf{R}_1\vert} \Vert}.
\end{equation}
We utilize a simple adaptive time step size strategy in this example
\begin{equation}
	\Delta t_{n+1} = \begin{cases}
      5 * 1e-2 & \text{if $\vert E(\Psi^{n+1}) - E(\Psi^{n})\vert \ge 1e-2$},\\
      5 * 1e-4 & \text{otherwise}.
    \end{cases}      
\end{equation}
The reason why we choose this strategy is to make full use of the unconditional energy stability. We employ a large time step size to let the solution close to the ground state as soon as possible, followed by a transition to a small time step size for accuracy. 
The nonuniform mesh size is defined by 
\begin{equation}
	h(x,y,z)= \begin{cases}
      0.2& \text{if $x^2+y^2+z^2 \le 1.8^2$},\\
      2.5 & \text{otherwise}.
    \end{cases}      
\end{equation}
And we set the thickness of the ball to 1.

\begin{figure}[h!]
	\subfloat{\includegraphics[scale = 0.22]{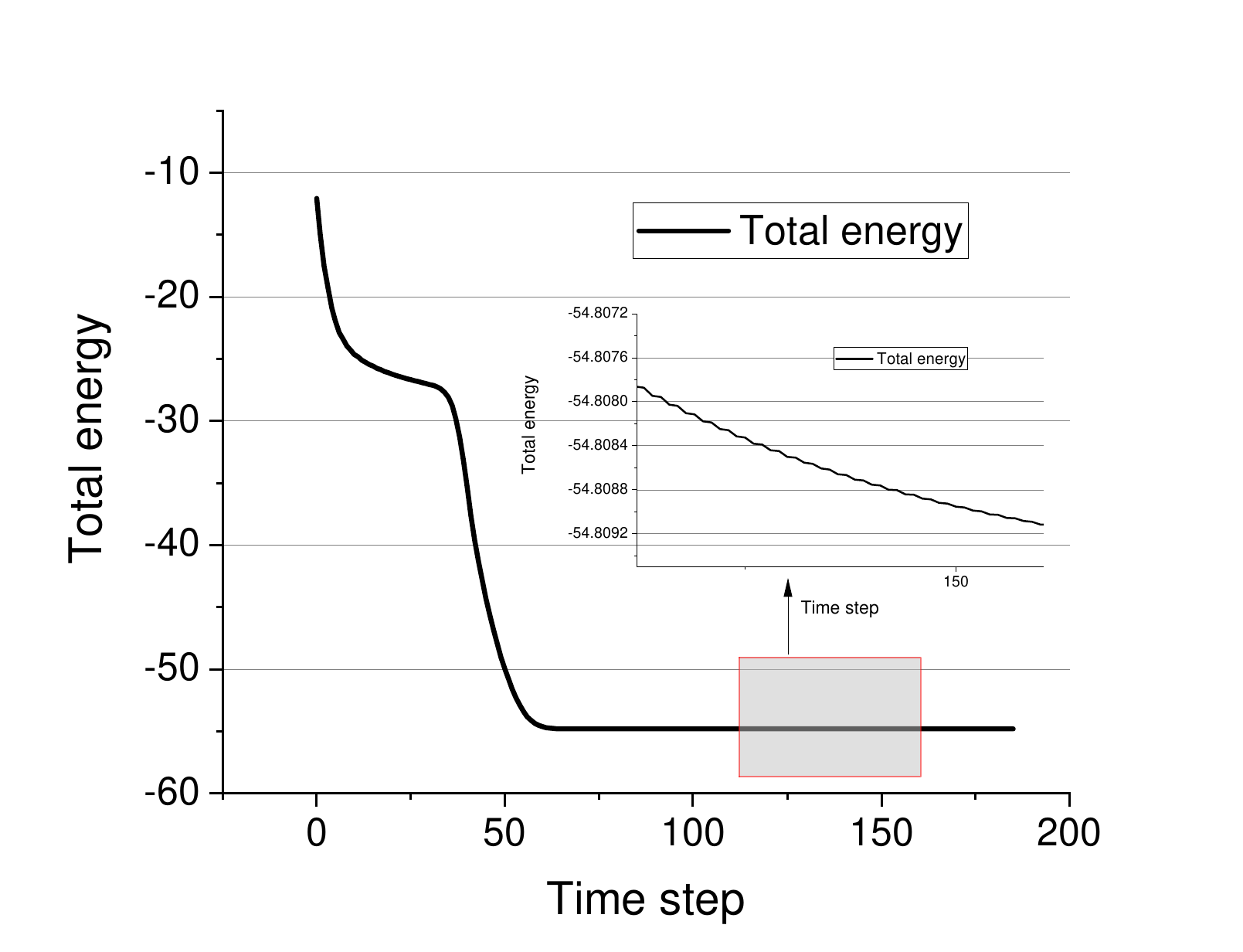}} 
	\hspace{0cm}
	\subfloat{\includegraphics[height = 5cm, width = 5cm]{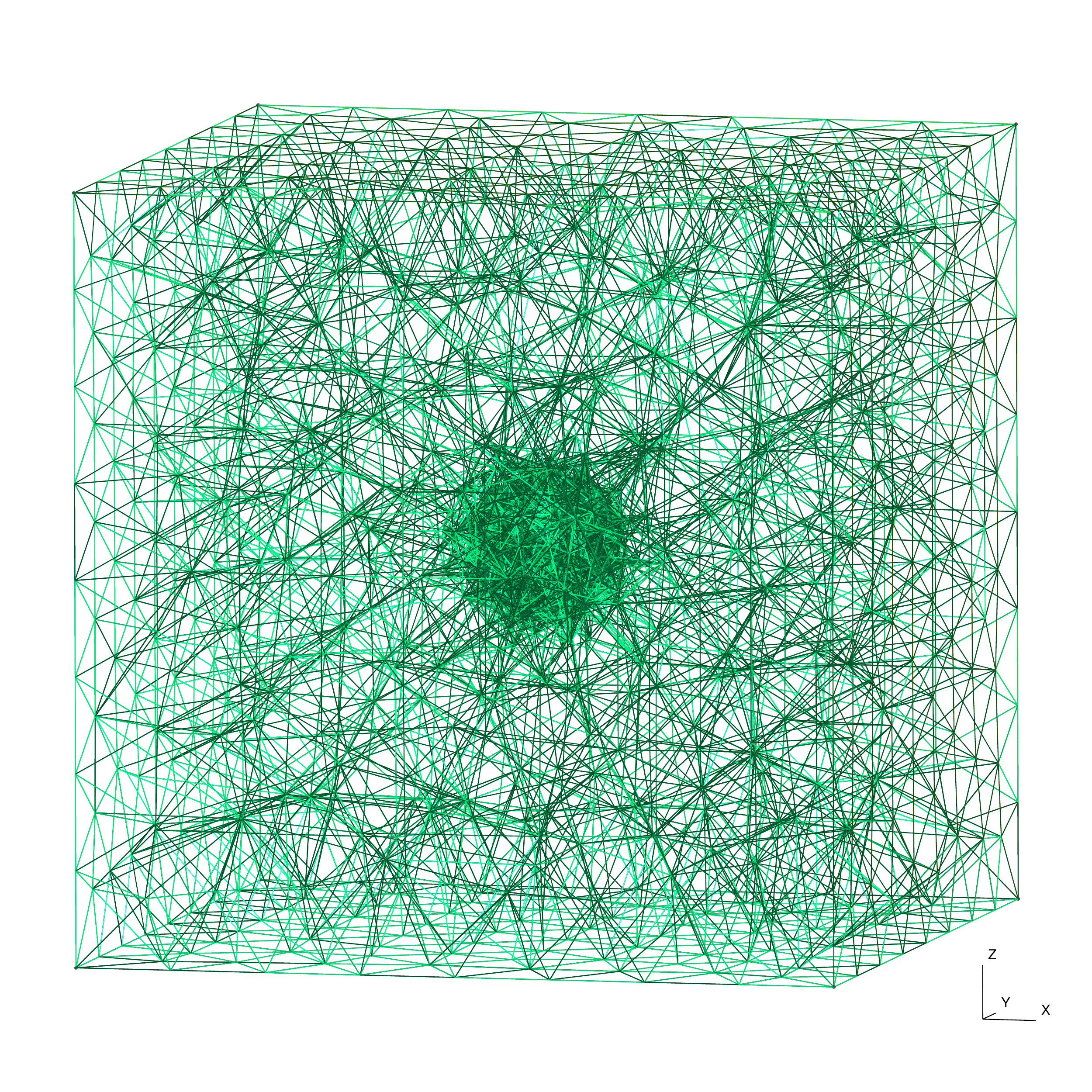}}
	\caption{The evolution of the computed total energy (in Hartree) with the time step (left) and the nonuniform mesh used for domain discretization with a total number of degree of freedoms 3323 (right) for the electronic structure of a methane molecule (Example \ref{Ex.3}). The evolution of the computed total energy demonstrates that our scheme is strictly energy stable even with large time steps. We purposely enlarge the energy scale at later time steps to show the total energy indeed decreases strictly, though slowly, after 100 time steps, verifying our theoretical analysis reported in this paper.}
	\label{Ex.3 Energy}
\end{figure}

\begin{figure}[h!]
	\centering
	\includegraphics[scale = 0.1]{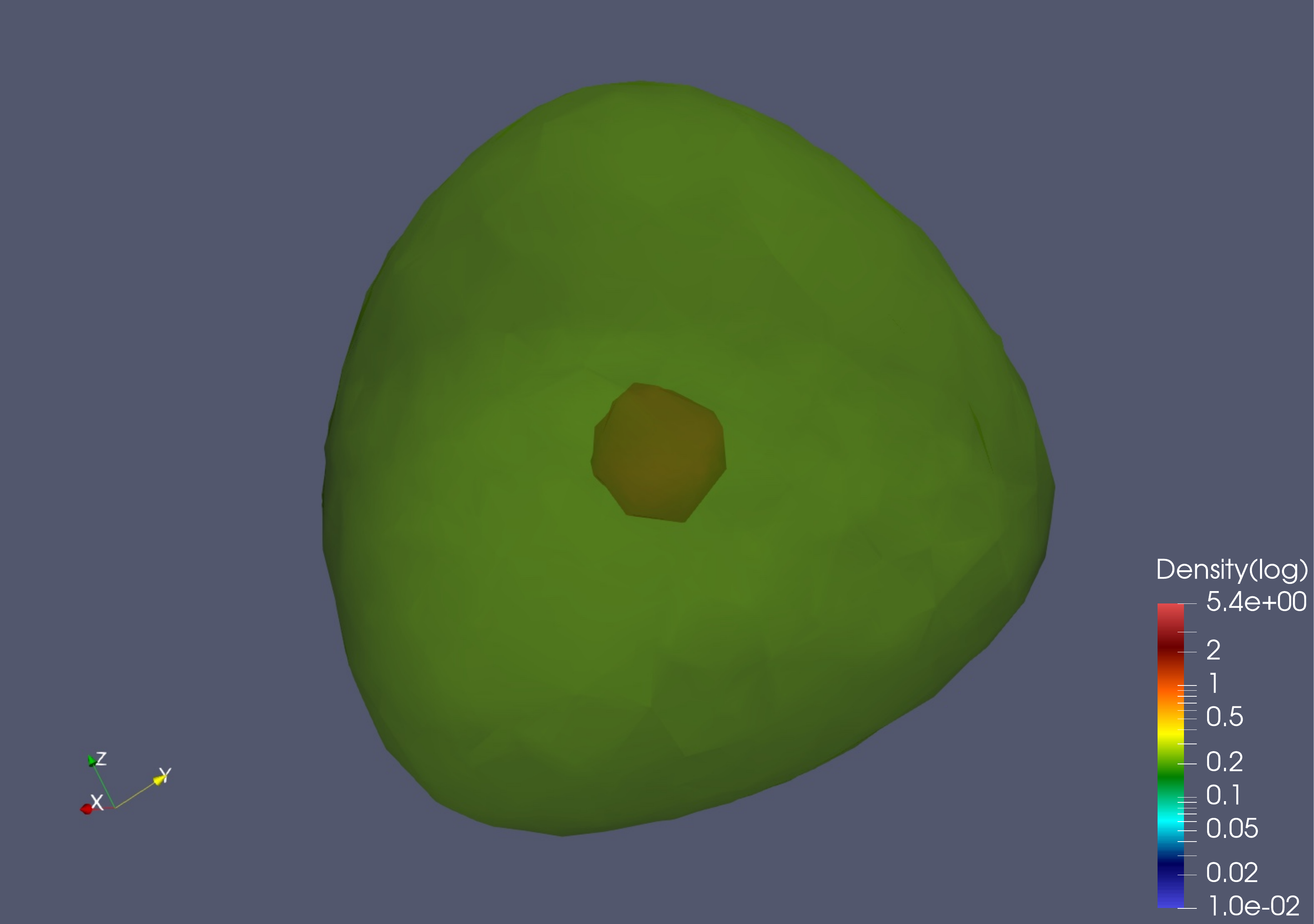}
	\caption{The 3D contour plot for the predicted electronic structure of a methane molecule (Example \ref{Ex.3}). The shown shape of the methane molecule is stretched by the four hydrogen atoms sitting outside. We purposely peek into the core to visualize the dense electronic cloud around the carbon atom sitting in the center.}
	\label{Ex.3 3D}
\end{figure}

\begin{figure}[h!]
	\subfloat{\includegraphics[scale = 0.1]{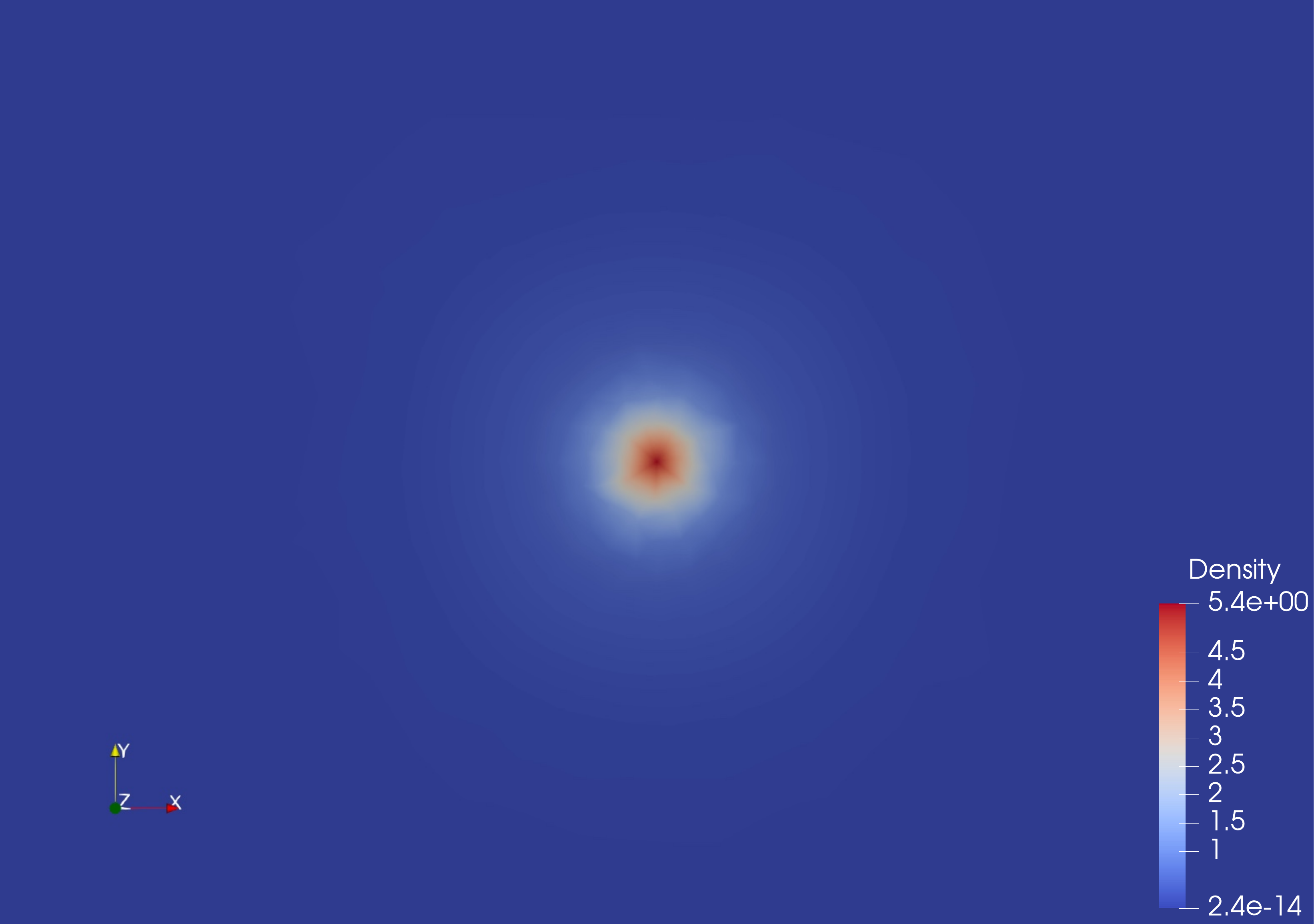}} 
	\hspace{1cm}
	\subfloat{\includegraphics[scale = 0.1]{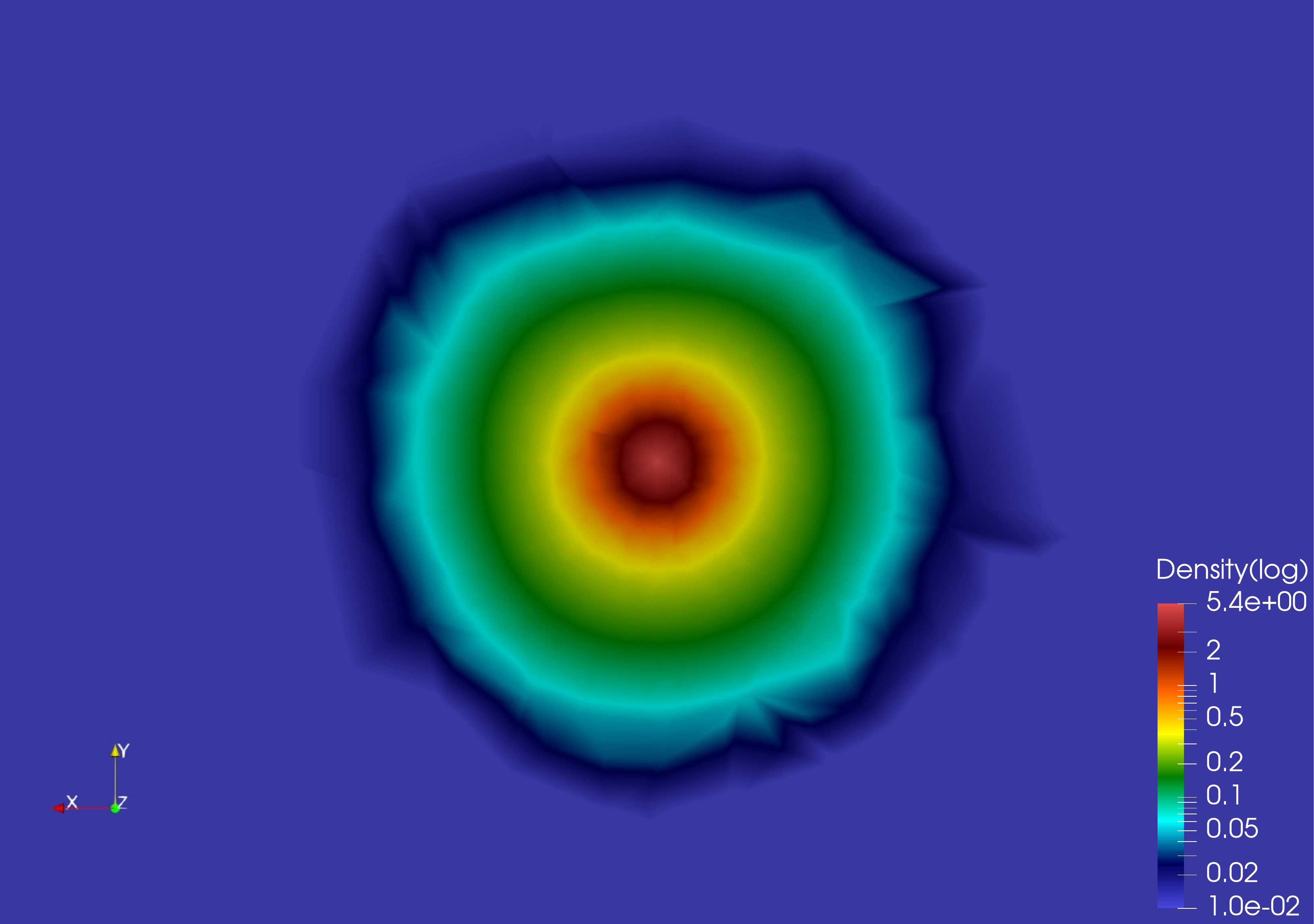}}
	\caption{Profiles of the computed electron density function in a linear scale (left) and in a log scale (right) for the $X Y$ cross-section of the electronic structure of a methane molecule (Example \ref{Ex.3}).}
	\label{Ex.3 result}
\end{figure}

\smallskip
We observe that the discrete energy dissipation law still holds from Figure \ref{Ex.3 Energy}. Analogously, the total energy travels a considerable distance from the initial condition to the ground state. 
We suggest the readers use the initial condition close enough to the ground state. The numerical examples are designed to show the power of unconditional energy stability. The 3D contour plot is presented in Figure \ref{Ex.3 3D}. We can see the dense electronic cloud around the carbon atom and the sparse cloud around the four hydrogen atoms (one at the top and three at the bottom).
The numerical density profiles are presented in Figure \ref{Ex.3 result}. 

\begin{remark}
	\normalfont We actually tested various meshes on the same problem, wherein it was discovered that the number of time steps remained almost the same across all tested meshes. This phenomenon shows that the proposed scheme inherits the property from the PDE, where the extended gradient has an exponential decay over time $t$ \cite{dai_gradient_2020}.
\end{remark}

\section{Conclusions}\label{Sect. 6}
This paper presents an unconditionally energy-stable, orthonormality-preserving, component-wise splitting iterative scheme for the Kohn-Sham gradient flow based model. The proposed scheme perfectly inherits the two desirable properties from the PDE model and offers an iterative approach that mitigates the need for significant computational resources. Rigorous proof and several numerical examples are given in the paper to verify, theoretically and numerically.  

Even though quite a number of orthonormality-preserving schemes have been proposed for the Kohn-Sham gradient flow based model, their energy stability typically depends on the time step size. To the best of our knowledge, the algorithm presented in this paper is the first unconditionally energy-stable algorithm reported in the literature that also preserves orthonormality exactly. The component-wise splitting technique simplifies the computation further and reduces the number of simultaneous equations to be solved at a time. The unconditional energy stability guarantees that we can use large time step sizes in the simulation, while component-wise splitting allows us to solve small systems in each time step, saving massive computational costs as well as memory requirements as compared to existing algorithms.  Unlike the SAV-based methods, our proposal method does not modify the original energy and does not have a tuned parameter for the algorithm; thus, our proposal method presents a robust and efficient choice of the Kohn-Sham gradient flow based solution procedures.  

\smallskip
The unconditionally energy-stable yet orthonormality-preserving iterative strategies presented in this paper shed light on the future study of the Kohn-Sham gradient flow based approach.  Future work can include applying adaptive strategies in mesh generation and time step sizes, developing proper preconditioners to accelerate the computation, energy-stable treatment of the exchange-correlation term, and many other improvements. The ultimate goal is to establish a linear version of the proposed scheme that still possesses the unconditional energy stability and the orthonormality-preserving property.  Many theoretical questions are also intriguing; for example, it is an interesting ongoing work to rigorously show the number of iterations in this unconditionally energy-stable yet orthonormality-preserving iterative scheme is independent of (or weakly depends on) the dimension of the Galerkin finite element spaces.

\vspace{0.15in}
\noindent
 {\small{\bf Acknowledgements}
This work is supported by King Abdullah University of Science and Technology (KAUST) through the grants BAS/1/1351-01 and URF/1/5028-01. The work of Huangxin Chen was supported by the NSF of China (Grant No. 12122115, 11771363).

\bibliographystyle{siam}
\bibliography{references}
\end{document}